\documentclass[sn-mathphys]{sn-jnl}

\usepackage{graphicx}%
\usepackage{multirow}%
\usepackage{amsmath,amssymb,amsfonts}%
\usepackage{amsthm}%
\usepackage{mathrsfs}%
\usepackage[title]{appendix}%
\usepackage{xcolor}%
\usepackage{textcomp}%
\usepackage{manyfoot}%
\usepackage{booktabs}%
\usepackage{algorithm}%
\usepackage[algo2e]{algorithm2e}

\usepackage{algorithmicx}%
\usepackage{algpseudocode}%
\usepackage{listings}%
\usepackage{flafter} 
\usepackage{enumitem}
\usepackage{amsfonts}
\usepackage{comment}
\RequirePackage[labelformat=simple]{subcaption}

\usepackage{blindtext}
\usepackage{cleveref}
\usepackage{multirow}
\usepackage{amssymb} 
\usepackage{minitoc}
\usepackage{lscape}
\usepackage{bbm, dsfont}

\theoremstyle{thmstyleone}%
\theoremstyle{thmstyletwo}%

\theoremstyle{thmstylethree}%

\begin{document}
\title[Article Title]{Transfer Learning in Bayesian Optimization for Aircraft Design}


\author*[1,2]{ \sur{Ali Tfaily}}\email{ali.tfaily@mail.mcgill.ca}


\author[5,6]{ \sur{Youssef Diouane}}\email{youssef.diouane@polymtl.ca}

\author[3,4]{ \sur{Nathalie Bartoli}}\email{nathalie.bartoli@onera.fr}

\author[7,1,6]{ \sur{Michael Kokkolaras}}\email{mkokkol@clemson.edu}

\affil*[1]{\orgdiv{Department of Mechanical Engineering}, \orgname{McGill University},  \country{Canada}}

\affil[2]{\orgdiv{Advanced Product Development}, \orgname{Bombardier}, \country{Canada}}

\affil[3]{\orgdiv{DTIS, ONERA}, \orgname{Université de Toulouse}, 
\country{France}}

\affil[4]{\orgdiv{Fédération ENAC ISAE-SUPAERO ONERA}, \orgname{Université de Toulouse}, 
\country{France}}

\affil[5]{\orgdiv{Department of Mathematical and Industrial Engineering}, \orgname{Polytechnique Montréal}, \orgaddress{ \country{Canada}}}

\affil[6]{\orgdiv{GERAD}, \orgaddress{ \country{Canada}}}

\affil[7]{\orgdiv{School of Mechanical and Automotive Engineering}, \orgname{Clemson University}, \orgaddress{ \country{USA}}}


\abstract{The use of transfer learning within Bayesian optimization addresses the disadvantages of the so-called \textit{cold start} problem by using source data to aid in the optimization of a target problem. 
We present a method that leverages an ensemble of surrogate models using transfer learning and integrates it in a constrained Bayesian optimization framework. We identify challenges particular to aircraft design optimization related to heterogeneous design variables and constraints. We propose the use of a partial-least-squares dimension reduction algorithm to address design space heterogeneity, and a \textit{meta} data surrogate selection method to address constraint heterogeneity. Numerical benchmark problems and an aircraft conceptual design optimization problem are used to demonstrate the proposed methods. Results show significant improvement in convergence in early optimization iterations compared to standard Bayesian optimization, with improved prediction accuracy for both objective and constraint surrogate models.}




\maketitle

\section{Introduction}
In the context of aircraft design, Bayesian optimization is a promising and growing field as presented in~\cite{bartoli2016improvement,roy2017mixed,lam2018advances,priem2020efficient,jim2021bayesian, saves2022bayesian}. However, Bayesian optimization suffers from slow local convergence in early iterations of some applications~\cite{priem2020efficient, bai2023transfer}.  One of the reasons for this slow convergence can be attributed to the so-called \textit{cold start} problem where each new Bayesian optimization is initiated without assuming any prior knowledge of the blackbox functions of the objective and constraints. 

Transfer learning can be beneficial for addressing the cold start issue in Bayesian optimization that contributes to slow convergence~\cite{niu2020decade, zhuang2020comprehensive, weiss2016survey, bai2023transfer}. 
Transfer learning in Bayesian optimization has been investigated and applied in the field of machine learning as surveyed by Bai \textit{et al.}~\cite{bai2023transfer} where the authors break down transfer learning in Bayesian optimization into four different components: (1) surrogate model design, (2) acquisition function definition, (3) initial design of experiments (DOE) strategy, and (4) design space selection. In this work, we focus on using transfer learning in the surrogate model design component of Bayesian optimization to address inaccurate surrogate model predictions in early iterations that slow down convergence considerably.  

Most Bayesian optimization applications use Gaussian processes (or variations thereof) as the surrogate model for the objective and constraint functions. In machine learning applications, Gaussian processes with a large number of design variables face computational and accuracy challenges. To that end, other surrogate modeling methods have been proposed, e.g., neural networks~\cite{springenberg2016bayesian}, deep Gaussian processes~\cite{wistuba2021few}, neural processes~\cite{wei2021meta}, and prior fitted networks~\cite{muller2023pfns4bo}. Such methods are not ideal for the types of problems considered here, where the number of design variables is less than 100 and the amount of data is relatively low. We will therefore focus on the use of Gaussian process surrogate models based on the methods presented in~\cite{tfaily2026transfer}.

The use of transfer learning on surrogate models of Bayesian optimization can be based on a single surrogate model that uses all the source and target data at once such as the approach proposed by~\cite{yogatama2014efficient, lindauer2018warmstarting}, or it can be based on a combination of surrogate models separated by source and target tasks as in~\cite{feurer2018practical, wistuba2018scalable}. Yogotama and Mann proposed a sequential model-based optimization that transfers information by constructing a common response surface for all source and target datasets~\cite{yogatama2014efficient}.
Feurer \textit{et al.} developed an ensemble model for Bayesian optimization that can incorporate the results of past optimization runs using so-called ranking-weighted Gaussian process ensembles (RGPEs)~\cite{feurer2018scalable}. They first train the Gaussian process posterior for every source task, and use the posterior mean function to compute the number of so-called misranked pairs between one source problem and the target problem. 
Wistuba \textit{et al.} proposed a two-stage transfer surrogate model (TST-R) for hyperparameter optimization of machine learning models~\cite{wistuba2016two}. The first stage of the surrogate model approximates the hyperparameter response functions of a new data set and each data set from the source data individually with Gaussian processes. The second stage combines the first-stage models by taking into account the similarity between the new target data and the source data from previous experiments. They construct a ranking of hyperparameter configurations as well as a prediction about the uncertainty of this ranking. They use the models from source problems to obtain the mean of a target problem. 
Another approach is the use of \textit{meta} data  that describe properties of the data used for the source and target problems~\cite{brazdil1994characterizing,rivolli2018characterizing}. Chen and Liu applied a weighted ensemble of Gaussian processes within Bayesian optimization (similar  to the one in~\cite{feurer2018practical}) on a one-dimensional turbomachinery problem~\cite{chen2022efficient}. 

Other research uses transfer learning to define a prior distribution of a Gaussian process using source data and then apply the prior to a target problem as proposed in~\cite{van2018hyperparameter} and~\cite{hellan2023data}. More recently, Mahboubi \textit{et al.} proposed a point-by-point transfer learning method that integrates mixtures of Gaussian experts with Bayesian optimization to address the cold start problem, demonstrating improvements on both synthetic and real-world datasets~\cite{mahboubi2025point}. Chen \textit{et al.} provided a comprehensive review of ensemble of surrogates methods in blackbox engineering optimization, highlighting their effectiveness in improving surrogate robustness but noting that most ensemble methods focus on combining different surrogate types rather than transferring knowledge across problems~\cite{chen2024ensemble}. In the area of constrained Bayesian optimization without transfer learning, Yu proposed a constraint-handling framework leveraging Prior-data Fitted Networks, a foundation transformer model that evaluates objectives and constraints simultaneously through in-context learning, demonstrating significant speedups over Gaussian process-based methods~\cite{yu2025constrained}.

Notably, existing transfer learning methods for Bayesian optimization focus exclusively on unconstrained problems or transfer only the objective function surrogate. None of the aforementioned approaches address the transfer of constraint surrogate models, which is critical for engineering design optimization where constraints often outnumber objectives and significantly impact convergence. Furthermore, these methods assume homogeneous design spaces between source and target problems, limiting their applicability to real-world scenarios where design variables and constraints may differ across optimization tasks. This work addresses these gaps by extending transfer learning to constrained Bayesian optimization with explicit handling of heterogeneous design spaces and constraints. 

This paper is organized as follows. In~\Cref{sec:TLBO-motivation}, we identify a set of motivating scenarios faced by aircraft designers when conducting aircraft design Bayesian optimization, and we present an aircraft conceptual design problem. We provide background on Bayesian optimization and transfer learning of surrogates in~\Cref{sec:BO-background}. \Cref{sec:TLBO-methodology} presents the method we propose for using transfer learning in Bayesian optimization, including the ensemble of surrogates method and novel approaches for handling heterogeneity. Numerical and aircraft design problem results are presented in~\Cref{sec:TLBO-results}, where we compare the optimization results of the proposed transfer learning algorithm with existing Bayesian optimization methods. Conclusions are drawn and perspectives are offered in~\Cref{sec:TLBO-conclusions}. 
\section{Aircraft design and transfer learning}
\label{sec:TLBO-motivation}
\subsection{Motivation and context}
Aircraft design Bayesian optimizations are initiated with a cold start, where designers trade  initial design of experiments (DOE) size (and amount of corresponding evaluations) with optimization convergence. A designer formulates an optimization problem by defining the objective function, constraint functions, and design variables. The typical Bayesian optimization (see Algorithm~\ref{EGOalgorithm} in~\Cref{sec:BO-background}) conducts a DOE, builds an initial surrogate model, and initiates the optimization loop. In scenarios where the optimization problem shares the same design variables, constraints, and objective function with an existing DOE, a designer may use that existing DOE to build the initial surrogate models and skip the cold start DOE step; this is the so-called \textit{hot start} Bayesian optimization~\cite{bartoli2019adaptive}. However, in most practical applications, designers would not repeat the same optimization problem with the same design space, output and their ranges. Therefore, such a hot start approach is typically not possible. 
For aircraft conceptual design optimization, we identify the following scenarios where a source optimization or DOE and target optimization differ such that a conventional hot start cannot be used: 
\begin{itemize}
    \item A change in aircraft design requirements such as aircraft range, cruise speed, balanced field length, approach speed, etc. A source optimization or DOE based on a set of design requirements can be different from the target optimization due to differing market requirements of the target optimization, management decisions based on the source optimization or DOE results, or new certification regulations that were not applicable to the source optimization or DOE.  
    \item A change in constraints driven by configuration or technological selections such as landing gear or flight control constraints. A source optimization may have been based on a  specific landing gear configuration (e.g., trailing arm type), a specific flight control system configuration (e.g., a slatted wing), or a specific flight control system actuation technology (e.g., linear electro-hydraulic servo-actuator); whereas a target optimization may have differing configurations or technologies (a cantilever landing gear type, an unslatted wing configuration, and a rotary electro-mechanical actuator respectively). Each configuration or technology may necessitate its own set of constraints which leads to a potential differences between the constraints and their limits between the source and target optimizations. 
    \item A change in fixed inputs, e.g., a fixed weight or cabin geometry. A target optimization may require a specific avionics function that is performed by a technology associated with a fixed input weight to the optimization, and this technology (and its associated weights) were not used in the source optimization setup. Similarly, a target optimization may change the cabin geometry relative to the source optimization which leads to several implications on an optimization's output (e.g., fuselage weight or tail moment arm). These different fixed inputs lead to different output values which prevents a designer from using a source optimization data in the target optimization. 
    \item A change of fidelity of models within the multi-disciplinary analysis (MDA). A source optimization may have been performed using a low fidelity aerodynamics model. Whereas a target optimization may use a newly calibrated low-fidelity aerodynamics model or a high-fidelity model. This change of model within the MDA leads to different results even for the same set of inputs between the source and target optimizations. 
    \item A change in simulation codes in the MDA between the source and target optimization. Updates to simulation codes due to bug fixes or changes in capabilities can lead to differing output for the same input. 
\end{itemize}
The motivation of this work is to enable transferring knowledge between source and target optimizations in the presence of changes in requirements, constraints, inputs, or models listed above. The aim is to reduce convergence time of new optimizations by creating transferable surrogate models of disciplines within an optimization.
Existing methods for transfer learning in Bayesian optimization do not modify surrogate models prior to transferring knowledge from the source to the target problem. In addition, previous work on the weighting of surrogate ensembles does not separate the design space prior to defining these surrogate ensembles, as explained in~\cite{tfaily2026transfer}. Existing methods for transfer learning using ensembles of Gaussian processes as surrogates do not address the need for accurate constraint models in Bayesian optimization where the value of predictions is just as important as the shape of the constraint functions. Finally, existing methods using ensembles of Gaussian processes as surrogates in Bayesian optimization do not handle heterogeneous design variables or heterogeneous constraints, i.e., when there is a mismatch between the design variables or the constraints output of the source and target problems.

The objective of this work is to develop algorithms to enable faster Bayesian optimization by applying transfer learning techniques while addressing the shortcomings of existing work presented above.
\subsection{An aircraft conceptual MDO problem}
\label{sec:into-aircraft-problem}
\begin{figure}[H]
    \centering
    \begin{subfigure}{0.45\textwidth}
    \centering 
    \hspace*{-0.68in}
    \includegraphics[width=7.5 cm]{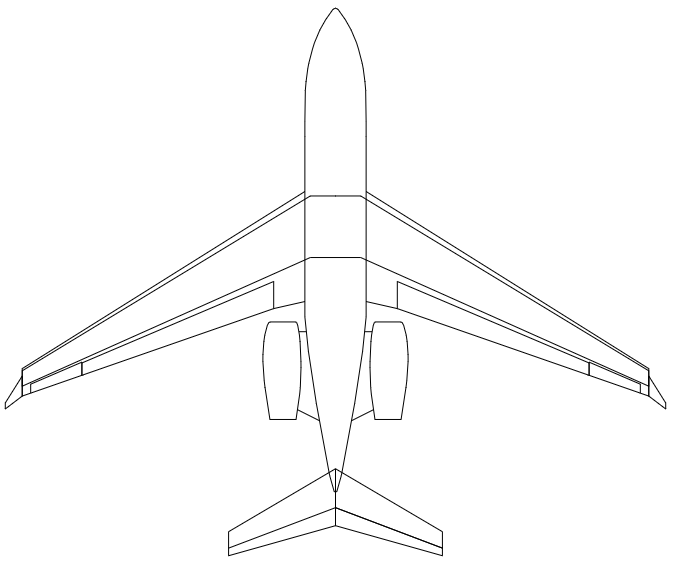}
    \caption{S2 aircraft}   
    \end{subfigure}
    \hfill
   \begin{subfigure}{0.45\textwidth}
    \centering 
    \includegraphics[width=7.5 cm]{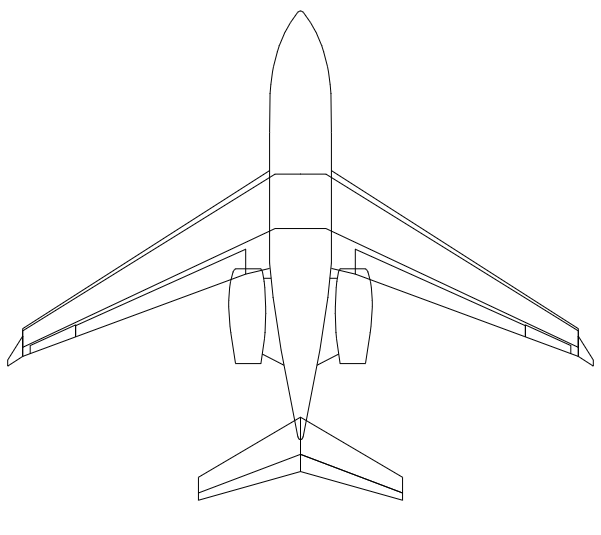}
    \caption{S3 aircraft}   
    \end{subfigure}
    \vfill
    \begin{subfigure}{0.45\textwidth}
    \centering 
    \hspace*{-0.3in}
    \includegraphics[width=7.5 cm]{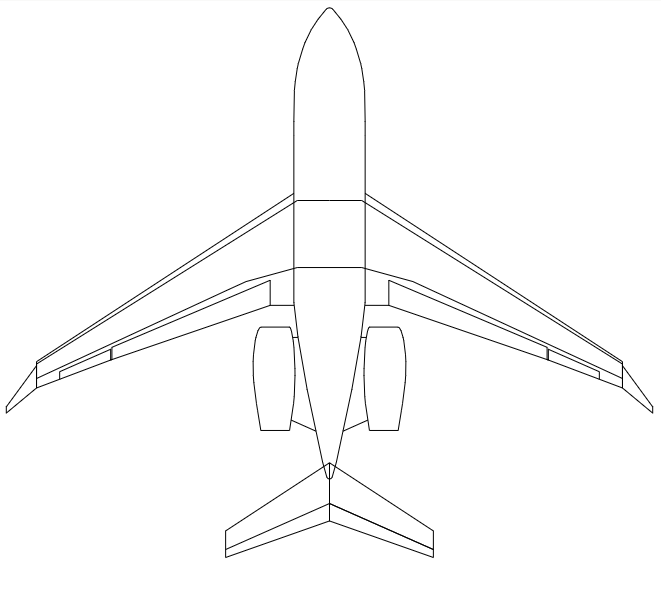}
    \caption{S4 aircraft}   
    \end{subfigure}
    \hfill
   \begin{subfigure}{0.45\textwidth}
    \centering 
    \includegraphics[width=7.5 cm]{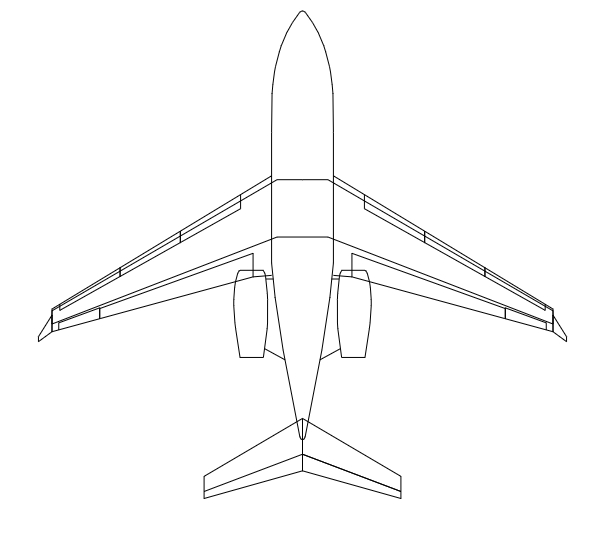}
    \caption{S5 aircraft}   
    \end{subfigure}
    \caption{Scaled top views of selected aircraft from the source problem's DOE highlighting the differences in wing and fuselage geometries.}
    \label{TLBO-ac-top-views-figure}
\end{figure}

We present an aircraft conceptual design optimization problem inspired by the examples identified previously. The problem of interest is a single target optimization problem. We consider available existing data based on five source DOEs (which are the results of previous optimizations). The definition of target and source problems is presented in~\Cref{TLBO-ac-design-prob-descriptions}  showing top views of selected aircraft from the DOE of each problem in~\Cref{TLBO-ac-top-views-figure}.

\begin{table}[h!]
    \caption{Setup of source and target problems in the aircraft conceptual MDO transfer learning problem.}
    \label{TLBO-ac-design-prob-descriptions}
    \begin{tabular}{|p{2cm}||p{8cm}|}
        \hline
         Problem & Description \\
         \hline \hline
         Target (T) & BRAC \\
         \hline
         Source 1 (S1)& BRAC with a modified engine weight (+500 lb) and a relaxed BFL constraint (+10\%) \\
         \hline
         Source 2 (S2)& BRAC with a longer fuselage (+100 inch) and longer range (+20\%)\\
         \hline
         Source 3 (S3)& BRAC with higher cruise speed requirement (+ Mach 0.05) \\
         \hline
         Source 4 (S4)& BRAC with a different engine architecture and rear spar configuration\\
         \hline    
         Source 5 (S5)& BRAC with slats \\
         \hline 
    \end{tabular}
\end{table}

\begin{figure}[H]
    \centering
	\includegraphics[scale = 0.35
 ]{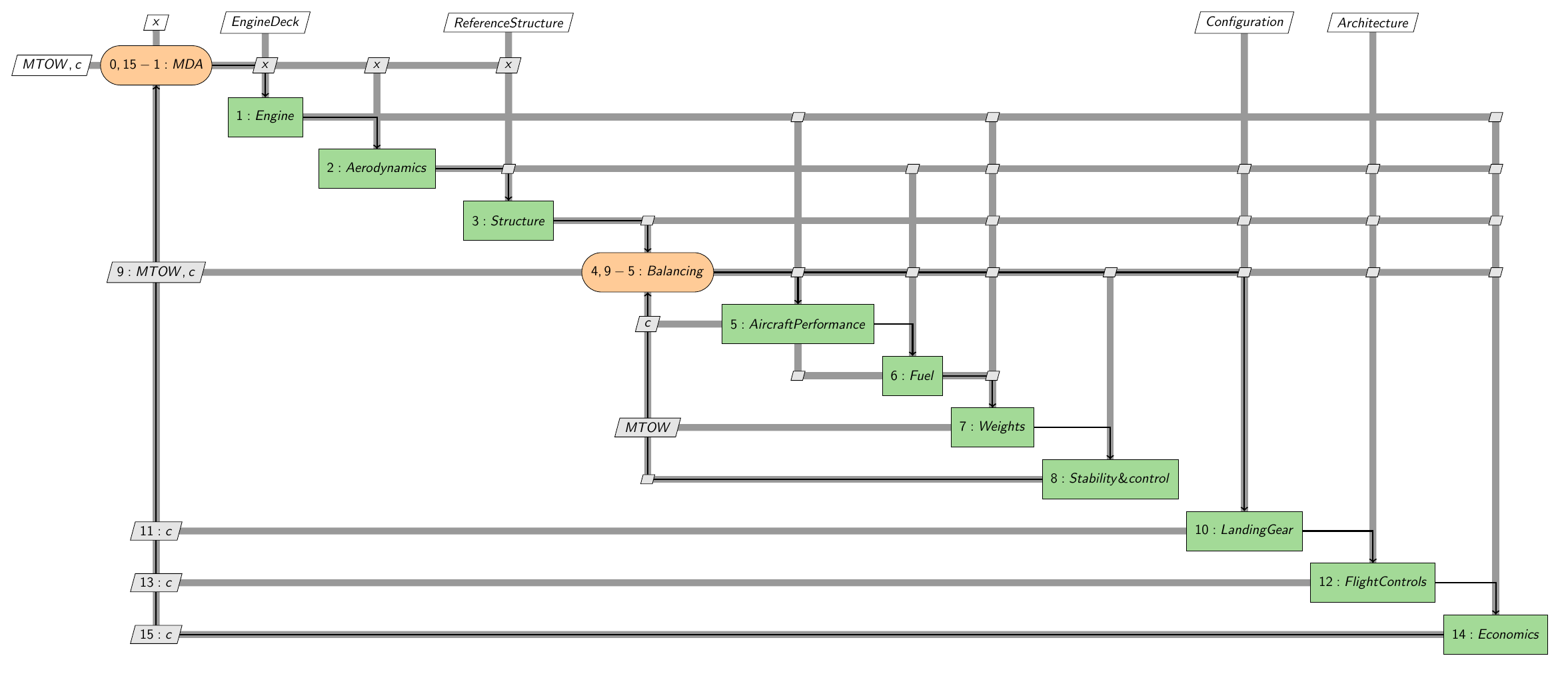}
    \caption{XDSM representation of the aircraft conceptual design MDA framework.}
    \label{xdsm}
\end{figure}

The target and source problems share a homogeneous objective function: minimize the maximum takeoff weight (MTOW) of the Bombardier Research Aircraft (BRAC)~\cite{reist2019cross, priem2020efficient} using heterogeneous requirements, constraints, and design variables.
The purpose is to leverage existing source data to improve the performance of the target optimization. 

\begin{equation*}\label{tlbo-cmdo-formulation-eq}
\Bigg\{
\begin{aligned}
& \underset{x \in \Omega}{\text{minimize}}
& & \text{MTOW}(x) \\
& \text{subject to}
& & c_i(x) \leq 0, \; i = 1, \ldots, m,
\end{aligned}
\end{equation*}
where $x\in \Omega \subseteq \mathbb{R}^{16}$ is the vector of bounded design variables listed in Table~\ref{TLBO-design-variables} and $c_i(x) \; i = 1, \ldots, 13$ are the inequality constraints described in~\Cref{TLBO-ac-constraints-definition}. We use Bombardier's conceptual multidisciplinary optimization (MDO) framework~\cite{piperni2013development} considering tube-and-wing aircraft configurations with a T-tail design and turbofan engines. The MDO framework includes an aircraft multidisciplinary analysis (MDA) environment which is comprised of sizing and simulation models of all major disciplines in aircraft conceptual design. The extended design structure matrix (XDSM)~\cite{lambe2012extensions} of the MDA environment is presented in Figure~\ref{xdsm}.
\begin{table}[h]
        \caption{ Design variables and their applicability (noted as a checkmark \checkmark) to the  target (T) and each source (S1 to S5).}
        \label{TLBO-design-variables}
        \begin{tabular}{|p{1.5cm}||p{4.7cm}||p{0.5cm}|p{0.5cm}|p{0.5cm}|p{0.5cm}|p{0.5cm}|p{0.5cm}|} 
             \hline
             Input  & Description & T & S1 & S2 & S3 & S4 & S5\\ [0.5ex] 
             \hline\hline
             $x_1$ & Rubber engine scaling factor & \checkmark & \checkmark & \checkmark$^1$ & \checkmark$^1$ & \checkmark & \checkmark \\ 
             \hline
             $x_2$ & Rubber engine bypass ratio &  &  &  &  & \checkmark &  \\ 
             \hline
             $x_3$ & Rubber engine overall pressure ratio &  &  &  &  & \checkmark &  \\ 
             \hline
             $x_4$ & Wing aspect ratio & \checkmark & \checkmark & \checkmark & \checkmark & \checkmark & \checkmark \\ 
             \hline
             $x_5$ & Wing area & \checkmark & \checkmark & \checkmark$^1$ & \checkmark$^1$ & \checkmark & \checkmark$^1$ \\ 
             \hline
             $x_6$ & Wing trailing edge sweep & \checkmark & \checkmark & \checkmark & \checkmark & \checkmark & \checkmark \\ 
             \hline
             $x_7$, $x_8$ & Wing rear spar chord-wise location & \checkmark & \checkmark & \checkmark & \checkmark & \checkmark & \checkmark \\  
             \hline
             $x_9$ & Wing sweep & \checkmark & \checkmark & \checkmark & \checkmark & \checkmark & \checkmark \\ 
             \hline
             $x_{10}$ & Wing taper ratio & \checkmark & \checkmark & \checkmark & \checkmark & \checkmark & \checkmark \\ 
             \hline
             $x_{11},\ldots,x_{14}$ & Wing thickness-to-chord ratios & \checkmark & \checkmark & \checkmark & \checkmark & \checkmark & \checkmark \\  
             \hline
             $x_{15}$ & Rear spar kink span-wise location &  &  &  &  & \checkmark &  \\ 
             \hline
             $x_{16}$ & Slat chord ratio &  &  &  &  &  & \checkmark \\ 
             \hline
        \end{tabular}
\footnotesize{$^1$ using different bounds for the same design variable}
\end{table}
\begin{table}[h]
        \caption{A list of the aircraft conceptual design problem constraints.}
        \label{TLBO-ac-constraints-definition}
        \begin{tabular}{|p{1.6cm}||p{4cm}||p{0.55cm}|p{0.55cm}|p{0.55cm}|p{0.55cm}|p{0.55cm}|p{0.55cm}|} 
             \hline
             Constraint  & Description & T & S1 & S2 & S3 & S4 & S5\\ [0.5ex] 
             \hline\hline
             $c_1(x)$ & Balanced field length (BFL) & = & \footnotesize{+10\%} & = & = &  = & = \\ 
             \hline
             $c_2(x)$ & Initial cruise altitude (ICA) & = & = & = & = & = & = \\ 
             \hline
              $c_3(x)$ & Aircraft reference speed ($V_{\mbox{ref}}$) & =  & = & = & = & = & = \\ 
             \hline
              $c_4(x)$ & Excess fuel weight & = & = & = & = & = & = \\ 
             \hline
              $c_5(x)$ &Inboard Wing flight controls actuation height clearance & $^*$ & $^*$ & $^*$ & $^*$ & $^*$ & $^*$ \\ 
             \hline
             $c_6(x)$ &Outboard Wing flight controls actuation height clearance & $^*$ & $^*$ & $^*$ & $^*$ & $^*$ & $^*$ \\ 
             \hline
            $c_7(x)$& Inboard Wing flight controls actuation chord clearance & $^*$ & $^*$ & $^*$ & $^*$ & $^*$ & $^*$ \\ 
            \hline
            $c_8(x)$ & OutboardWing flight controls actuation chord clearance & $^*$ & $^*$ & $^*$ & $^*$ & $^*$ & $^*$ \\ 
             \hline
              $c_9(x)$& Wing chord clearance for landing gear integration & $^*$ & $^*$ & $^*$ & $^*$ & $^*$ & $^*$ \\ 
             \hline
              $c_{10}(x)$ & Wing tip chord & = & = & = & = & = & = \\ 
             \hline
             $c_{11}(x)$& Aircraft climb performance & = & = & = & = & = & = \\ 
             \hline
             $c_{12}(x)$ & Aircraft mission range & = & = & \footnotesize{+20\%} & = & = & = \\  
             \hline
             $c_{13}(x)$ & Slat actuator height & \footnotesize{N/A} & \footnotesize{N/A} & \footnotesize{N/A} & \footnotesize{N/A} & \footnotesize{N/A} & $^*$ \\  
             \hline
        \end{tabular}
\footnotesize{$^*$ Constraints defined using internal simulation models within the MDA}\\
\footnotesize{$=$ Equivalent constraint bounds for target and source problems}\\
\footnotesize{$\%$ Heterogeneous constraint thresholds that vary by a fixed percentage relative to the target problem}
\end{table}

Note that the reference engine weight and fuselage length are fixed inputs to the MDA; therefore, they do not appear in the list of design variables in~\Cref{TLBO-design-variables}. In addition, the cruise speed requirement is set in the MDA and the aircraft performance simulation model analyzes its impact on the mission, i.e., cruise speed is not a constraint.

\section{Bayesian optimization and transfer learning}
\label{sec:BO-background}
Bayesian optimization is thought of having its origins in Efficient Global Optimization (EGO)~\cite{jones1998efficient}. It is a surrogate-assisted optimization approach that has gained significant attention in the last two decades. It targets problems with expensive or noisy blackbox functions, promising reduced computational cost and efficient convergence using a probabilistic approach. Consider the unconstrained optimization problem
\begin{equation}
\underset{x \in \Omega}{\text{min}} \quad y(x),
\label{eq:unconstrained-opt}
\end{equation}
where $y(x)$ is a noisy blackbox function and $x$ is its input vector within a design space $ \Omega \subseteq \mathbb{R}^{D}$. In Bayesian optimization, a surrogate $\hat{y}(x)$ is created based on the input and output of the function $y(x)$ and is used in a sequential optimization framework. There are three components in Bayesian optimization: (1) an initial DOE, (2) a surrogate model, and (3) an acquisition function that is used to select the points to be evaluated using the blackbox function $y(x)$.

\subsection{Gaussian process surrogate models}
In Bayesian optimization, the most popular surrogate models are Gaussian processes~\cite{williams2006gaussian}, also known as Kriging~\cite{krige1951statistical}. Let us assume a DOE with $N$ samples evaluated using a blackbox function $y(x)$ and denote its input vectors of dimensions $d$ as $X_N = \{x_n\}_{n=1}^{N}$ and its output vector as $Y_N = \{y_n\}_{n=1}^{N}$ with $x_n \in \mathbb{R}^D$ and $y_n \in \mathbb{R}$, $\forall n \in N$. A Gaussian process is a stochastic process that approximates $y(x)$ at any point $x$ by a Gaussian random variable $\hat{Y}(x)$ which is defined by its mean and covariance functions, $\hat{y}(x)$ and $k(\cdot,\cdot)$ respectively. The mean function is defined as
\begin{equation}
\hat{y}(x) = \hat{\beta} + r^T(x)K^{-1}(Y_N - \mathbf{1}\hat{\beta})
\label{eq:GP-mean}
\end{equation}
and the variance is defined as
\begin{equation}
\hat{s}^2(x) = \hat{\sigma}^2(x)(1 - r^T(x)K^{-1}r(x)),
\label{eq:GP-variance}
\end{equation}
where $\hat{\beta}$ is a mean bias, $r(x) = \{k(x, x_1), \ldots, k(x, x_N)\}$, $K$ is an $N \times N$ covariance matrix of the combination of components of the DOE $K_{ij} = k(x_i, x_j)$, $\mathbf{1}$ is an $N \times 1$ column vector of 1's, and $\hat{\sigma}$ is the variance of the covariance function. In this work, we use the squared exponential covariance function:
\begin{equation}
k(x_i, x_j) = \hat{\sigma}^2 \prod_{d=1}^{D} \exp\left(\frac{(x_i^d - x_j^d)^2}{2l_d}\right),
\label{eq:SE-kernel}
\end{equation}
where $x_i$ and $x_j \in \mathbb{R}^D$, the length scale vector $l = \{l_1, \ldots, l_D\} \in \mathbb{R}$, the variance $\hat{\sigma}$, and the mean bias $\hat{\beta}$ are obtained by the training data of the DOE. The most common method to obtain these hyperparameters is the maximum likelihood estimation (log-ML)~\cite{forrester2008engineering}:
\begin{equation}
\text{log-ML} = -\frac{1}{2}\left(N\ln(2\pi\hat{\sigma}^2) + \ln(\det K) + \frac{(Y_N - \mathbf{1}\hat{\beta})^T K^{-1}(Y_N - \mathbf{1}\hat{\beta})}{\hat{\sigma}^2}\right).
\label{eq:log-ML}
\end{equation}

\subsection{Acquisition functions}
The acquisition function is a critical part of the Bayesian optimization framework. Many acquisition functions have been proposed and investigated in literature. Expected Improvement (EI)~\cite{jones1998efficient} is one of the most widely used:
\begin{equation}
\text{EI}(x) = \int_{-\infty}^{\infty} I(x)\phi(z) \, dz,
\label{eq:EI}
\end{equation}
where $I(x)$ is the improvement at point $x$ relative to the the minimum value of the objective function observed so far $y_{\min} -  \hat{y}(x)$, \(y_{\min}\) is the minimum value of the objective function observed so far, \({\hat{y}(x)}\) and \({\hat{s}(x)}\) are mean and standard deviation of the  Gaussian process, and $\phi$ is the probability distribution function, and $z = (y_{\min} - \hat{y}(x))/\hat{s}(x)$.

\subsection{Constraint handling}
\label{sec:constraint-handling}
We consider an inequality constrained optimization formulation. We use the Bayesian optimization extension to handle constrained blackbox optimization problems as proposed by~\cite{schonlau1997computer} where they propose creating surrogate models of the inequality constraints $c_1(x), \ldots, c_m(x)$ similar to the objective function and apply these constraints on the acquisition function maximization problem. We define the constrained optimization problem as
\begin{equation}
\begin{aligned}
& \underset{x \in \Omega}{\text{min}} && y(x) \\
& \text{subject to} && c_i(x) \leq 0, \quad i = 1, \ldots, m.
\end{aligned}
\label{eq:constrained-opt}
\end{equation}
The BO problem is then defined as
\begin{equation}
\begin{aligned}
& \underset{x \in \Omega}{\text{max}} && \alpha(x) \\
& \text{subject to} && \hat{c}_i(x) \leq 0, \quad i = 1, \ldots, m,
\end{aligned}
\label{eq:acq-constrained}
\end{equation}
where $\alpha(x)$ is the selected acquisition function, and $\hat{c}_1(x), \ldots, \hat{c}_m(x)$ are surrogate models of the inequality constraints.
The overall algorithm for a constrained Bayesian optimization framework is defined in Algorithm~\ref{EGOalgorithm}.

\begin{algorithm}[!h]
\SetAlgoLined
\LinesNumbered
\setcounter{AlgoLine}{0}
\KwIn{Initial DOE}
\KwOut{Best feasible point from the DOE}
\For{\(i=1, \ldots, \mbox{max\_iter}\)}{
    Build Gaussian process surrogate models of the objective and constraint functions.\;
    
    Maximize an acquisition function to find \(x^{i+1}\).\;
    
    Evaluate objective and constraint functions at \(x^{i+1}\).\;
    
    Update the DOE.\;
 } 
 \caption{Constrained Bayesian optimization framework.}
 \label{EGOalgorithm}
\end{algorithm} 
\subsection{Ensemble of surrogates using transfer learning}
\label{sec:ensemble-TL}

We use an ensemble of Gaussian process surrogate models in lieu of the conventional Gaussian processes used in Bayesian optimization based on the work presented in~\cite{tfaily2026transfer}. The method starts by transferring the posterior of source models to target models as follows. Then, surrogate models are weighted based on multiple criteria according to the intended function. 
\subsubsection{Weighted ensemble of surrogates}
A source model ($\hat{y}_{\mbox{s}}^n$, $\hat{s}_{\mbox{s}}^n$) is adjusted by considering both a scaling factor and a bias, as presented in~\cite{tfaily2026transfer} 
\begin{equation}
\begin{aligned}
& \hat{y}_{\mbox{s}'}^n(x) = \alpha_n \hat{y}_{\mbox{s}}^n(x)+ \beta_n, \\
\end{aligned}
\label{TL-prior-alphab}
\end{equation}
where the scaling term $\alpha_n$ and the bias term $\beta_n$ are constant parameters that can be fitted to the target data $\mbox{DOE}_{\mbox{T}}$ to create a modified set of source models using a prespecified  error metric  $\varepsilon$ (e.g., mean squared error). The corresponding adjusted variance is then calculated as follows
\begin{equation}
\begin{aligned}
& \hat{s}_{\mbox{s}'}^n(x) = \alpha_n \hat{s}_{\mbox{s}}^n(x). \
\end{aligned}
\label{TL-prior-alphab-s}
\end{equation}

All source models go through this adjustment step as presented in~\Cref{{TL-prior-transfer-algo}}. 
\begin{algorithm} [h!]
\SetAlgoLined
\LinesNumbered
\KwIn{$N$ Source datasets $\mbox{DOE}_{\mbox{S}}^{1,...,N}$ and target dataset $\mbox{DOE}_{\mbox{T}}:(\mbox{X}_{\mbox{T}}, \mbox{Y}_{\mbox{T}})$.}
\KwOut{ {Modified source GP models ($\hat{y}_{\mbox{s}'}^{1,...,N}$, $\hat{s}_{\mbox{s}'}^{1,...,N}$).}}

\For{\(n=1, \ldots, N\)}{
   Build source GP model ($\hat{y}_{\mbox{s}}^{n}$, $\hat{s}_{\mbox{s}}^n$) based on $\mbox{DOE}_{\mbox{S}}^{n}$ using~\Cref{eq:GP-mean,eq:GP-variance,eq:SE-kernel,eq:log-ML}.\;
    
    Use the posterior of the $\hat{y}_{\mbox{s}}^{n}$ to calculate $\hat{y}_{\mbox{s}}^n (\mbox{X}_{\mbox{T}})$.\;

    Calculate $(\alpha_n,  \beta_n) = \underset{\alpha,  \beta \in \mathbb{R}}{\mbox{argmin }} \varepsilon (\alpha,\beta)$ where $\varepsilon$ is an error metric  between \hbox{$\hat{y}_{\mbox{s}'} (\mbox{X}_{\mbox{T}}) = \alpha \hat{y}_{\mbox{s}}^n(\mbox{X}_{\mbox{T}})+ \beta$ and  $\mbox{Y}_{\mbox{T}}$} for a given $\alpha,  \beta \in \mathbb{R}$ \;

    {Compute the hyperparameters $\sigma_{0}^n$ and $l^n$.\;

     Set the modified source  GP model:  $\hat{y}_{\mbox{s}'}^n= \alpha_n \hat{y}_{\mbox{s}}^n + \beta_n $ and variance $\hat{s}_{\mbox{s}'}^n = \alpha_n \hat{s}_{\mbox{s}}^n $.}\;
 }
 
 \caption{Source surrogate model with posterior information~\cite[Algorithm 1]{tfaily2026transfer}.}
 \label{TL-prior-transfer-algo}
\end{algorithm} 
The method then combines the modified source surrogate models into an ensemble using a weighted combination approach. The ensemble prediction $\hat{y}$ at a point $x$ is given by:
\begin{equation}
\hat{y}(x) = \sum_{j=1}^{N} P(\hat{y}_j = \hat{y}^*) \hat{y}_j(x),
\label{eq:ensemble-prediction}
\end{equation}

where $N$ is the number of source models, $\hat{y}^*$ represents the true underlying function, and $P(\hat{y}_j = \hat{y}^*)$ is the probability that source model $j$ is the correct model. The computation of these probabilities is described in~\Cref{sec:multi-criteria}. This formulation was written for a mixture of experts approach with multiple clusters in~\cite{tfaily2026transfer} where the design space is divided into $K$ clusters and model probabilities are computed per cluster; however, in this work we use a single cluster ($K=1$) to minimize computational costs.

To obtain the standard deviation $\hat{s}$ of the ensemble of surrogates $\hat{y}$, we present two methods:
\begin{enumerate}
    \item Target surrogate model variance $\hat{s}_{\mbox{t}}$ calculated using~\Cref{eq:GP-variance}:
    \begin{equation}
    \hat{s}^2(x) = \hat{s}_{\mbox{t}}^2(x).
    \label{eq:variance-target}
    \end{equation}
    \item Weighted ensemble of surrogates variances:
    \begin{equation}
    \hat{s}^2(x) = \sum_{j=1}^{N} (P(\hat{y}_j = \hat{y}^*))^2 \hat{s}_j^2(x),
    \label{eq:variance-ensemble}
    \end{equation}
    where $\hat{s}_j(x)$ is the standard deviation of the $j$th source model.
\end{enumerate}

\subsubsection{Multi-criteria weighted ensembles}
\label{sec:multi-criteria}
To compute the model probabilities $P(\hat{y}_j = \hat{y}^*)$ in Equation~\eqref{eq:ensemble-prediction}, a multi-criteria ranking approach was used. The selectable criteria are based on the shape, accuracy, and variance of the surrogates in an ensemble.
A summation of weighted multi-criteria $C(\hat{y}_j = \hat{y}^*)$ to calculate a score for a model $\hat{y}_j$ that being the correct model is presented as:
\begin{equation}
C(\hat{y}_j = \hat{y}^*) = \sum_{l=1}^{N_C} w_l c_l(\hat{y}_j = \hat{y}^*),
\label{eq:multi-criteria-score}
\end{equation}
where $w_l \in [0, 1]$ for $l = 1, \ldots, N_C$ are the weights assigned to the selected criteria such that $\sum_{l=1}^{N_C} w_l = 1$, $c_l$ is the measure of the selected criterion, and $N_C$ is the number of measured criteria.

For the shape criterion, the discordant pairs method is used from~\cite{feurer2018practical,tfaily2026transfer}. The Epanechnikov quadratic kernel~\cite{epanechnikov1969non} is used as follows:
\begin{equation}
\kappa_\rho(\chi_S, \mbox{DOE}_{\mbox{T}}) = \delta\left(\frac{\tau(\chi_S, \mbox{DOE}_{\mbox{T}})}{\rho}\right)
\label{eq:epanechnikov}
\end{equation}
with
\begin{equation}
\delta(t) = \begin{cases}
\frac{3}{4}(1 - t^2) & \text{if } t \leq 1 \\
0 & \text{otherwise}
\end{cases}
\label{eq:delta}
\end{equation}
and
\begin{equation}
\tau(\chi_S, \mbox{DOE}_{\mbox{T}}) = \frac{\sum_{m=1}^{n}\sum_{k=2}^{n} \mathbf{1}((y_m^s < y_k^s) \oplus (y_m^t < y_k^t))}{n},
\label{eq:tau}
\end{equation}
where $\rho > 0$ is a predefined bandwidth, $\mbox{DOE}_{\mbox{T}}$ is the input-output data of the target problem $T$, and $\chi_S$ are the data predictions $\mbox{Y}_{\mbox{S}}$ on the inputs $\mbox{X}_{\mbox{T}}$ of $\mbox{DOE}_{\mbox{T}}$ using a source model $\mbox{S}$ with $n$ points.

For the accuracy criterion, a measure of the absolute relative error metric is used:
\begin{equation}
\tau_a(\chi_S, \mbox{DOE}_{\mbox{T}}) = \frac{1}{n}\sum_{m=1}^{n} \mathbf{1}\left(\frac{|y_m^s - y_m^t|}{|y_m^t|} > \epsilon_{\max}\right),
\label{eq:tau-accuracy}
\end{equation}
where $\epsilon_{\max}$ is a predefined maximum relative error.
Similarly, for the variance criterion:
\begin{equation}
\tau_v(\chi_S, \mbox{DOE}_{\mbox{T}}) = \frac{1}{n}\sum_{m=1}^{n} \mathbf{1}\left(\frac{\sigma_m^s}{y_{\max}^t} > \sigma_{\max}\right),
\label{eq:tau-variance}
\end{equation}
where $\sigma_{\max}$ is a predefined maximum relative variance, $y_{\max}^t$ is the maximum observed absolute output value in $\mbox{DOE}_{\mbox{T}}$, and $\sigma_m^s$ is the variance of point $m$ measured using source model S.

The probability of a source model $j$ being the correct model is then estimated as:
\begin{equation}
P(\hat{y}_j = \hat{y}^*) \approx \frac{C(\hat{y}_j = \hat{y}^*)}{\sum_{n=1}^{N} C(\hat{y}_n = \hat{y}^*)}
\label{eq:probability-model}
\end{equation}
where $N$ is the number of source models. 

The choice of criteria weights depends on the role of the surrogate model in the optimization. For the objective function surrogate, shape similarity is prioritized because the acquisition function relies on the relative ranking of predictions to identify promising regions; accurate absolute values are less critical since the optimizer seeks improvement directions rather than exact function values. For constraint surrogates, accuracy becomes equally important because constraint satisfaction depends on the actual predicted values relative to constraint bounds---a surrogate that correctly ranks points but systematically over- or under-predicts may misclassify feasibility. The variance criterion is included in both cases to favor source models with lower uncertainty, improving confidence in the ensemble predictions.

\subsubsection{Ensemble of surrogates algorithm}
The ensemble of surrogates using the transfer learning algorithm can be summarized in Algorithm~\ref{TL-ensemble-algo}. This algorithm is a simplified instance of ~\cite[Algorithm 3]{tfaily2026transfer} where only one cluster is regarded. 

\begin{algorithm}[!h]
\SetAlgoLined
\LinesNumbered
\setcounter{AlgoLine}{0}
\KwIn{Source datasets $\mbox{DOE}_{\mbox{S}}^{1,\ldots,N}$, target dataset $\mbox{DOE}_{\mbox{T}}$, selected criteria $N_C$, and criteria weights $w_l$ and bandwidths $\rho_l$ for $l \in 1, \ldots, N_C$}
\KwOut{Ensemble of surrogates using transfer learning}
Build target surrogate model $\hat{y}_t$ using $\mbox{DOE}_{\mbox{T}}$\;

\For{$n = 1, \ldots, N$}{
    Build $n$th source surrogate model $\hat{y}_{s_n}$ using $\mbox{DOE}_{\mbox{S}}^n$\;
    
    Use the posterior of $\hat{y}_{s_n}$ and $\mbox{DOE}_{\mbox{T}}$ to create a modified source surrogate model $\hat{y}_{s_n'}$ using Algorithm~\ref{TL-prior-transfer-algo}\;
    
    Calculate the criteria scores $C(\hat{y}_{s_n'} = \hat{y}^*)$ using Equation~\eqref{eq:multi-criteria-score}\;
    
    Calculate the probability $P(\hat{y}_{s_n'} = \hat{y}^*)$ using Equation~\eqref{eq:probability-model}\;
}
Return the ensemble of surrogates model $(\hat{y}, \hat{s})$ using Equations~\eqref{eq:ensemble-prediction} and \eqref{eq:variance-target} or \eqref{eq:variance-ensemble}\;
\caption{Ensemble of surrogate models using transfer learning.}
\label{TL-ensemble-algo}
\end{algorithm}

\section{Methodology}
\label{sec:TLBO-methodology}
In this section, we describe the algorithms and methods we propose for heterogeneous transfer learning in Bayesian optimization. Section~\ref{sec:heterogeneous-design-space} proposes an approach to handle heterogeneous design spaces using partial-least-squares dimension reduction. Section~\ref{sec:heterogeneous-constraints} presents a method to handle heterogeneous constraints based on meta data, which are features that describe the design variables and constraints, e.g., variable names. In Section~\ref{sec:ensemble-TL}, we present the ensemble of surrogates using transfer learning approach, including the transfer of prior method, weighted ensemble formulation, and multi-criteria weighting scheme. Section~\ref{sec:TLBO-algorithm} presents the overall algorithm to perform transfer learning in constrained Bayesian optimization.

\subsection{Heterogeneous design space}
\label{sec:heterogeneous-design-space}
\begin{figure}[!h]
    \centering
	\includegraphics[width=12 cm]{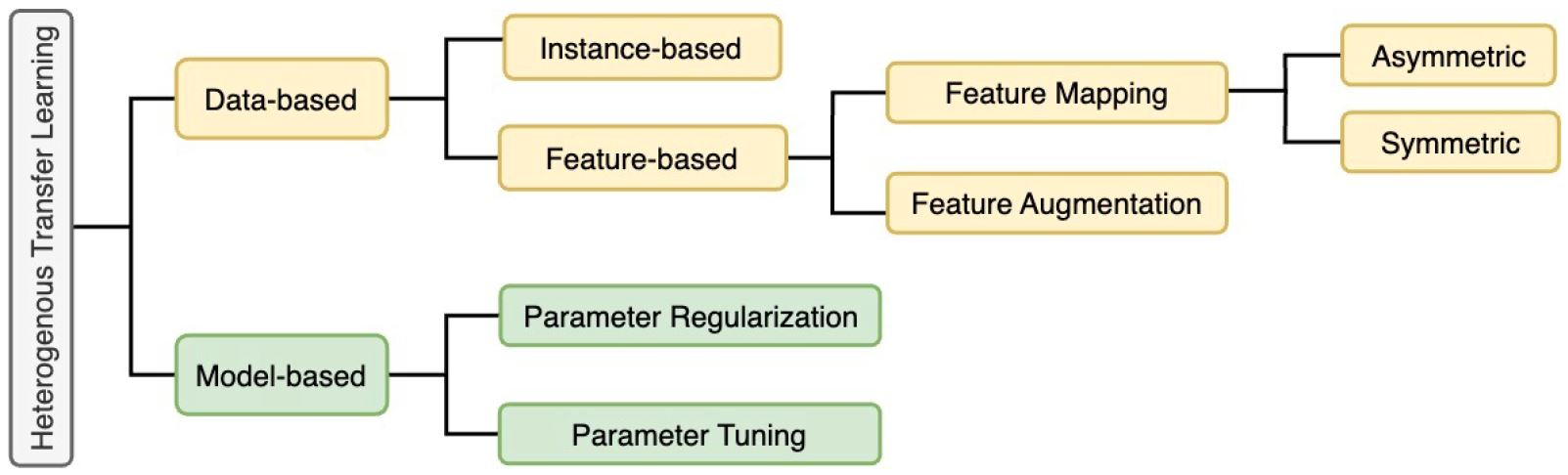}
    \caption{Heterogeneous transfer learning categories (adapted from~\cite{bao2023survey}).}
    \label{TL-hetero-categories-figure}
\end{figure}

Heterogeneity in transfer learning is encountered when the design space or the outputs of the source model(s) do not match the design space or outputs of the target model~\cite{weiss2016survey,zhuang2020comprehensive, bao2023survey, khan2024heterogeneous}. In real-world scenarios, the use of existing data typically comes with a variation of the design variables used between source and target problems. For example, in aircraft design settings source data may contain an engine by-pass ratio as a design variable whereas a target problem may have a fixed engine size and architecture. In this case, the target problem still has the engine by-pass ratio as a fixed parameter that can be used to solve the heterogeneity issue. However, in certain cases, source data may include design variables that do not even exist as fixed parameters in the target problem. For example, source data may be obtained from an aircraft that has slats with a design variable for the ratio of slats chord to the wing chord. Assuming that the target problem is an unslatted aircraft, this design variable does not exist and therefore cannot be used as a fixed parameter. In~\cite{bao2023survey}, Bao \textit{et al.} surveyed methods that handle heterogeneity in transfer learning problems and categorize the types of algorithms as shown in Figure~\ref{TL-hetero-categories-figure}. They distinguished between data-based and model-based methods. Data-based methods involve the transfer 
of either the original data or their transformed design variables to a target design space, allowing the target models to be trained with this augmented data, thereby enriching the available data within the target design space. 

Conversely, model-based methods center around constructing models within the source design space then applying target data to adapt model results to the target design space. In literature, heterogeneous transfer learning in Bayesian optimization is presented in~\cite{min2020generalizing} and~\cite{fantransfer}. In~\cite{min2020generalizing}, Min \textit{et al.} proposed a model-based approach that builds a neural network to estimate  design variables mapping from a source problem to a target problem using a neural network. In~\cite{fantransfer}, Fan, Han, and Wang proposed a neural network model that estimates a prior of the target model based on the source models and design variables. The neural network is trained using source design variables and posteriors of source surrogate models. In this work, the use of data-based methods is more appropriate for ensembles of surrogates since we are interested in maintaining each surrogate model as a separate entity. 

\begin{figure}[!h]
    \centering
        \begin{subfigure}{\textwidth}
            \centering           
            \includegraphics[height=3cm]{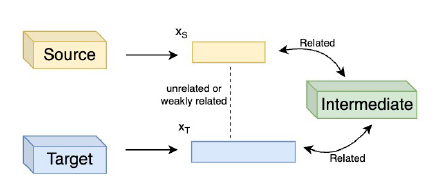}
            \caption{Instance-based method.}
            \label{HTL-instance-based-fig}            
        \end{subfigure}
        \vfill
        \begin{subfigure}{\textwidth}
        \centering
            \includegraphics[height=3.3cm]{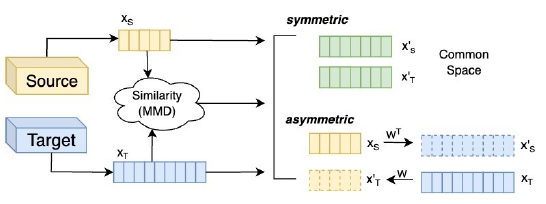}
            \caption{Feature mapping method.}
            \label{HTL-feature-map-fig}
    \end{subfigure}  
     \vfill
    \begin{subfigure}{\textwidth}
    \centering
            \includegraphics[height=3cm]{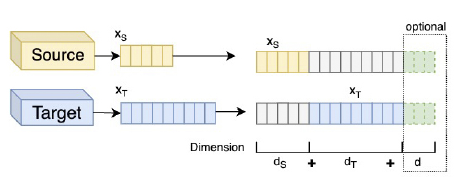}
            \caption{Feature augmentation method.}
            \label{HTL-feature-aug-fig}
    \end{subfigure}   
     \caption{Heterogeneous data-based transfer learning methods (adapted from \cite{bao2023survey}).}
     \label{HTL-data-based-fig}
\end{figure}

Assuming a set of design variables of a source problem $\mbox{X}_{\mbox{S}}$ and a set of design variables of a target problem $\mbox{X}_{\mbox{T}}$, Bao \textit{et al.} categorized the differences between the types of data-based transfer learning methods, i.e., instance or feature based, in Figure~\ref{HTL-data-based-fig}~\cite{bao2023survey}. Instance-based methods are best suited for situations where $\mbox{X}_{\mbox{S}}$ and $\mbox{X}_{\mbox{T}}$ are weakly related. An intermediate set of design variables is defined to build both the source and target surrogate models. The concept of feature-based methods is to increase the similarity between the source and target design variables. A common measure of similarity is the so-called Maximum Mean Discrepancy (MMD) that assesses the distances between the means of distributions between the source and target models~\cite{gretton2006kernel}. Feature mapping methods entail modifying the set of design variables of the source and target problems into a common design space, so-called symmetric, or projected design space on the source or target problems, so-called asymmetric. Feature augmentation methods add dimensions $d$ to the dimensions of the source and target problems, $d_{\mbox{S}}$ and $d_{\mbox{T}}$ respectively~\cite{li2013learning}.
In this work, we propose the use of a partial-least-squared dimension reduction algorithm to address heterogeneous design spaces between source and target problems similar to feature mapping methods in~\cite{bao2023survey} since the source and target problems will share a number of design variables in our aircraft optimization settings as we describe in~\Cref{TLBO-ac-design-prob-descriptions}. In addition, we are interested in reducing the computational cost of the surrogate models in Bayesian optimization. Bouhlel \textit{et al.} proposed a method for dimension reduction of Gaussian processes using partial least squares~\cite{bouhlel2016improving}. They developed a kernel that reduces the number of parameters of a Gaussian process. The key assumption is that some input variables have a less significant impact on the output of a Gaussian process than others. Therefore, they proposed to project the design space parameters onto smaller dimensions while monitoring the correlation between input variables and the output. The partial least squares method builds a relationship between these input variables and the projected variables using linear combinations of the input variables. An input vector \mbox{X} of $d$ dimensions is used to create $t$ principal components with $h$ dimensions to predict the output of a function $y$. These principal components represent the new coordinate system and are defined for $l=1,...,h$ as
\begin{equation}
    t^l = \mbox{X}\mbox{w}_*^{(l)},
\end{equation}
where $\mbox{w}_*^{(l)}$ for $l=1,...,h$ represent the weights of the input vector $\mbox{X}$ on principal component $t$. The vector $\mbox{w}_*^{(l)}$ is obtained by maximizing the covariance of $\mbox{X}^{(l-1)^t}\mbox{y}^{(l-1)}\mbox{y}^{(l-1)^t}\mbox{X}^{(l-1)}$. The number $h$ of principal components is chosen by minimizing the leave-one-out cross-validation error. Consequently, a point $x$ is transformed into principal component $t$ as follows
\begin{equation}
    x \longmapsto \bigl[\mbox{w}_{*1}^{(l)}x_1,...,\mbox{w}_{*d}^{(l)}x_d  \bigr]^t.
\end{equation}
Finally, the authors obtain the following Gaussian exponential kernel between points $x$ and $x'$:
\begin{equation}
    \begin{split}
    k(x,x')=&  \sigma ^ {2}\prod_{l=1}^{h}\prod_{i=1}^{d} \bigl[-\theta_{l}(\mbox{w}_{*i}^{(l)}x_i-\mbox{w}_{*i}^{(l)}x'_i)^2\bigr],
    \end{split}
\end{equation}
where $\sigma^2$ and $\theta _ {l} \in  [0,+  \infty[$ are the covariance and hyper-parameters of the Gaussian, respectively. Bouhlel \textit{et al.}  applied their method on numerical problems obtaining significant computational gain while maintaining prediction accuracy compared to that of conventional Gaussian processes. 

\begin{figure}[!h]
    \centering
	\includegraphics[width=12 cm]{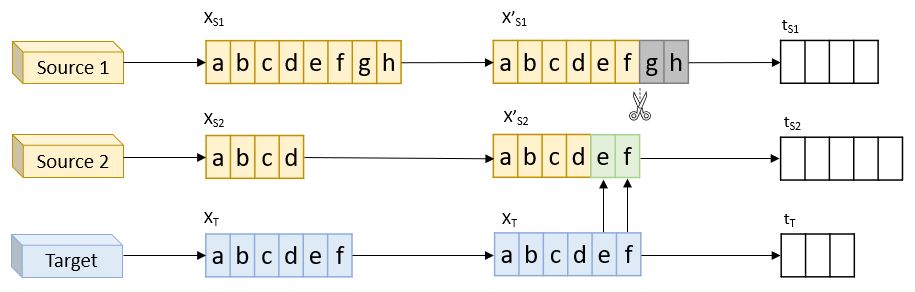}
    \caption{Heterogeneous design space in transfer learning using meta data.}
    \label{TL-hetero-PLS-figure}
\end{figure}

In the aircraft design optimization problems considered here, source problems share a number of design variables with the target problem. Input and output variables are also defined in advance. Therefore, we can assign meta data that define each input variable as presented in Figure~\ref{TL-hetero-PLS-figure}, where a transfer learning problem with two source input vectors and a target input vector is considered. Knowing these meta data, we propose to modify the KPLS algorithm in~\cite{bouhlel2016improving} by adding or removing variables to the modified input vectors of source problems to match the input variables of the target problem. Consequently, when performing the optimization to obtain the weights vector of the principal components $\mbox{w}_*^{(l)}$ for $l=1,...,h$, we assign a NULL value for  $\mbox{w}_*^{(l)}$ of each added or removed input variable of the source problems. The illustration in Figure~\ref{TL-hetero-PLS-figure} shows that input vectors $\mbox{X}_{\mbox{S1}}, \mbox{X}_{\mbox{S2}}, \mbox{X}_{\mbox{T}}$ are assigned defining meta data (a,b,c...). Therefore, when creating the KPLS components, variables with meta data g and h, $\mbox{w}_{\mbox{g}}^{(l)}$ and $\mbox{w}_{\mbox{h}}^{(l)}$ are assigned NULL values when creating the modified input vector $\mbox{X'}_{\mbox{S1}}$ of source problem 1 for all $l=1,...,h$. Similarly, $\mbox{w}_{\mbox{e}}^{(l)}$ and $\mbox{w}_{\mbox{f}}^{(l)}$ are assigned NULL values when creating the modified input vector $\mbox{X'}_{\mbox{S2}}$ of source problem 2. The resulting principal components $t_{\mbox{S1}}, t_{\mbox{S2}}, t_{\mbox{T}}$ will depend on the relationships between each modified input vector and the corresponding output. We note that for this approach, the target problem does not need to apply KPLS dimension reduction for predictions as opposed to the source problems.

Another approach to solve the heterogeneous design space problem in transfer learning of Bayesian optimization is the use of Hierarchical Gaussian Processes~\cite{saves2022general}. The idea is to combine all source and target data into a single Gaussian process and to add a set of so-called dimensional variables that affect the dimensions of a problem and decide if other so-called decreed variables, i.e. conditionally active or inactive variables, are acting or non-acting~\cite{audet2023general}. We determine that the use of hierarchical Gaussian processes is not ideal for solving the heterogeneity problem of transfer learning in Bayesian optimization since the need to add all source and target data points into a single Gaussian process significantly increases computational training needs at every iteration compared to the ensemble of surrogates approach where multiple source Gaussian processes are prebuilt and only the target Gaussian process is reconstructed at every iteration.

\subsection{Heterogeneous constraints}
\label{sec:heterogeneous-constraints}
\begin{figure}[!h]
    \centering
\includegraphics[width=6cm]{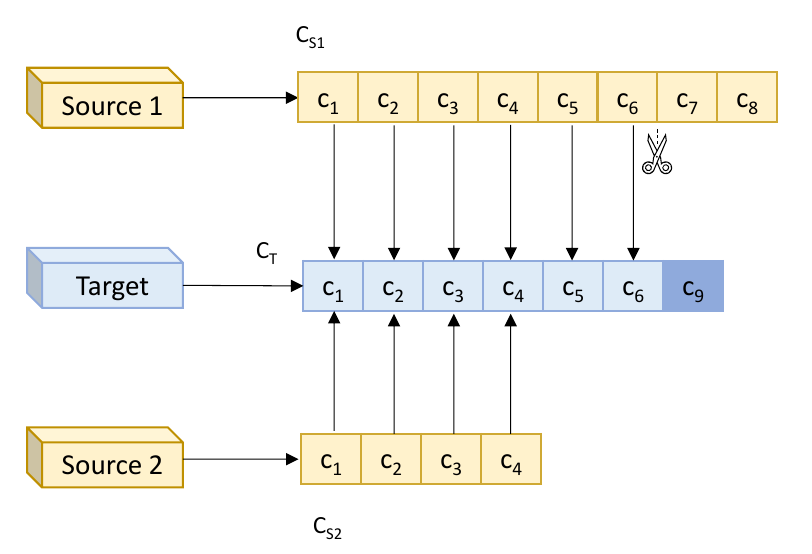}    
    \caption{Heterogeneous constraints transfer of surrogates approach showing the matching of source and target constraints using meta data.}
    \label{TL-hetero-constraints-figure}    
\end{figure}

In~\Cref{sec:TLBO-motivation}, we considered scenarios where constraints of source problems are not equivalent to target problem constraints. In this section, we present an approach to handle heterogeneous optimization constraints for transfer learning. Let's first consider a target problem with constraints vector $\mbox{C}_{\mbox{T}}$ and two source problems with constraints vectors $\mbox{C}_{\mbox{S1}}$ and $\mbox{C}_{\mbox{S2}}$ as shown in~\Cref{TL-hetero-constraints-figure}. One could use the ensemble of surrogates approach defined for the objective function using all the constraints data from both of the source problems to define each of the target constraint surrogates. However, for scenarios where there is a large number of constraints and a large number of source problems, such an approach may not be practical due to the additional computational cost needed to create ensembles for each target problem constraint. Therefore, we propose the use of meta data to select the source constraints surrogate models to be used for the ensemble of surrogates for a target constraint.

Assuming there exist meta data to describe each constraint, e.g., parameter name, we simply propose to match the same parameters from the source problems when creating the target problem constraint ensemble of surrogates. Therefore, in the schematic shown in~\Cref{TL-hetero-constraints-figure}, target constraints $c_1-c_4$ are built using source constraint models from both $\mbox{C}_{\mbox{S1}}$ and $\mbox{C}_{\mbox{S2}}$ whereas constraints $c_5$ and $c_6$ are built solely based on source constraint models from $\mbox{C}_{\mbox{S1}}$. 

\begin{figure}[!h]
    \begin{subfigure}{0.45\textwidth}
    \centering           
    \includegraphics[width=5 cm]{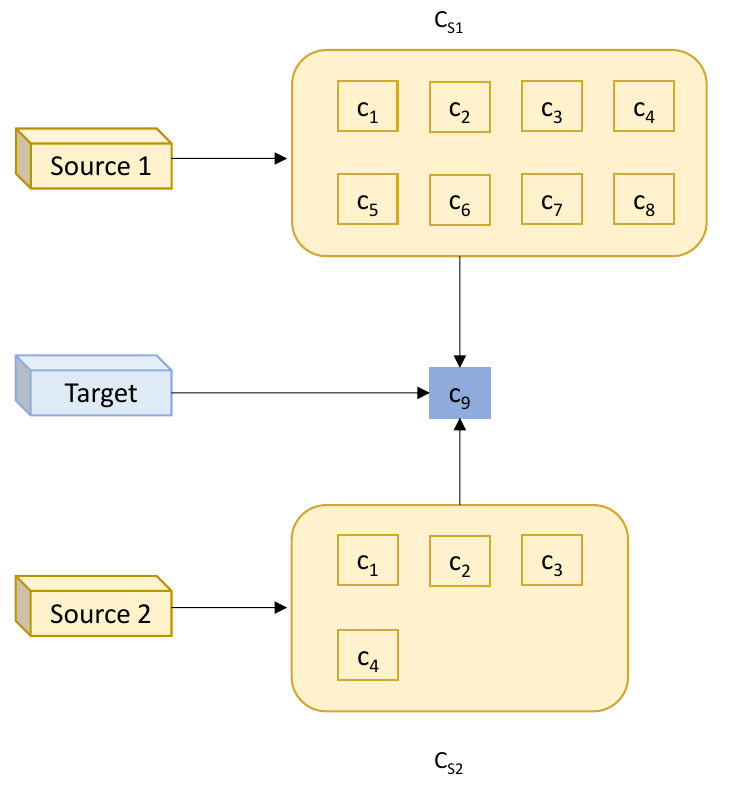}
    \caption{Use of all source constraint models for the ensemble of surrogates.}
    \label{TL-hetero-new-constraints-general-figure}      
    \end{subfigure}
        \hfill
    \begin{subfigure}{0.45\textwidth}
    \centering
	\includegraphics[width=5 cm]{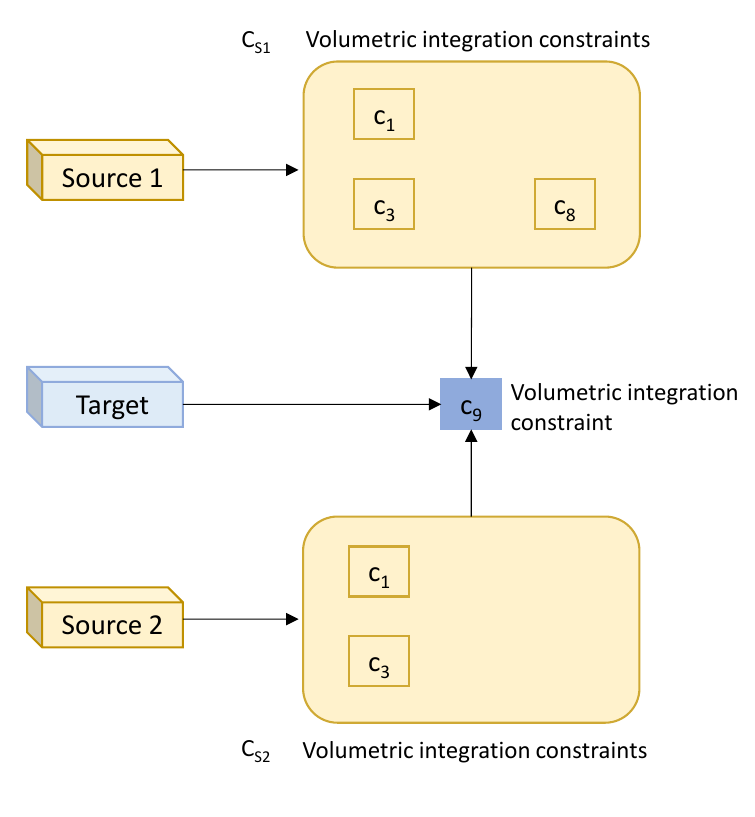}
    \caption{Use of subset of source  constraint models for the ensemble of surrogates.}
    \label{TL-hetero-new-constraints-meta-figure}       
    \end{subfigure}
\label{TL-hetero-new-constraints-figure}
\caption{Target problem constraints without matching meta data in source problems.}
\end{figure}

We also present a case where the target constraint does not exist in any of the source constraint data in~\Cref{TL-hetero-constraints-figure}, i.e., $c_9$ does not have a matching parameter in $\mbox{C}_{\mbox{S1}}$ nor $\mbox{C}_{\mbox{S2}}$. In such a case, we can use all the source problems constraints for the ensemble of surrogates of $c_9$ as presented in~\Cref{TL-hetero-new-constraints-general-figure} which may be computationally costly as explained above. Alternatively, we propose creating categories for constraints to reduce the number of models generated for the ensemble, as presented in~\Cref{TL-hetero-new-constraints-meta-figure}. For aircraft conceptual design problems, we categorize optimization constraints by discipline, as per the examples presented in~\Cref{table-aircraft-MDO-constraints1}. Let's assume that $c_9$ is a volumetric integration constraint in aircraft conceptual design, therefore, only parameters under the same category can be used to generate the ensemble of surrogates for $c_9$ as we present in~\Cref{TL-hetero-new-constraints-meta-figure}.

\begin{table}
    \caption{Aircraft MDO constraints.}
    \label{table-aircraft-MDO-constraints1}
    \begin{tabular}{ |p{2.2cm}|p{10cm}| }     
    \hline
    
    Constraints type & Constraints description\\
    \hline
    \hline
    \multirow{3}{4em}{Performance} & Minimum range at different speeds and payloads  \\
    & Maximum takeoff (or balanced) field length   \\
    & Maximum landing field length  \\ 
    & Minimum initial cruise altitude  \\ 
    & Maximum approach (or reference) speed \\
    & Maximum climb duration  \\     
    \hline
    \multirow{4}{5em}{Volumetric integration} & Minimum flight control actuation space envelope  \\ 
    & Minimum landing gear integration space envelope  \\
    & Minimum empty fuselage volume (outside of payload volume)  \\ 
    & Minimum fuselage underfloor height for avionics integration  \\
    & Maximum wing span  \\ 
    & Minimum excess fuel volume after design mission \\
    \hline
    \multirow{3}{4em}{Operational} & Maximum external noise in different flight phases (flyover, approach, takeoff)  \\ 
    & Maximum operating cost  \\ 
    & Ability of engines to provide power for anti-icing \\ 
    \hline    
    \multirow{1}{4em}{Environmental} & Emissions , noise, and life cycle analysis  \\ 
    \hline
    \end{tabular}
\end{table}


\subsection{Transfer learning for constrained Bayesian optimization}
\label{sec:TLBO-algorithm}
The overall algorithm for a constrained Bayesian optimization framework is defined in~\Cref{TLBO-algo}. Variations of this algorithm can be created to cater for the problems at hand. For example, building the probabilities for the ensemble of surrogates can be created prior to the Bayesian optimization loop to reduce the computational cost of the algorithm (or fixing these probabilities after a set number of iterations). In this work, we aim to maximize the use of target data, so we elect to re-calculate the  model probabilities in the ensemble of surrogates after every new evaluation point of the target functions (objective and constraints). 
\begin{algorithm}[!h]
\SetAlgoLined
\LinesNumbered
\setcounter{AlgoLine}{0}
\KwIn{Source datasets $\mbox{DOE}_{\mbox{S}}^{1,...,N}$ and initial target dataset $\mbox{DOE}_{\mbox{T}}$ }
\KwOut{Best feasible point from the target DOE}
Build source Gaussian processes for the objective and constraint functions $\hat{y}_{s}^{1,...,N}$ using $\mbox{DOE}_{\mbox{S}}^{1,...,N}$

\For{\(i=1, \ldots, \mbox{max\_iter}\)}{
    Build a transfer learning ensemble of Gaussian processes using~\Cref{TL-ensemble-algo} for the objective function and each constraint\;

    Maximize an acquisition function to find \(x^{i+1}\). \;
   
    Evaluate objective and constraint functions of the target model at \(x^{i+1}\).\;
    
    Update $\mbox{DOE}_{\mbox{T}}$.\;
 } 
 \caption{Constrained Bayesian optimization with transfer learning.}
 \label{TLBO-algo}
\end{algorithm}

\section{Numerical results}
\label{sec:TLBO-results}
We present analytical and aircraft conceptual design optimization  examples to demonstrate the proposed methods. 
\subsection{Implementation details}
The code is developed in Python 3 building on the following tool boxes: Scikit-learn v1.5.2~\cite{scikit-learn}, SMT 2.0~\cite{SMT2ArXiv}, SMAC3~\cite{lindauer-jmlr22a}, RGPE~\cite{feurer2018practical}, TST-R~\cite{wistuba2016two}, and the Bayesian optimization library in~\cite{gardner2014bayesian}. All results are obtained using an Intel® Xeon® CPU E5-1650 v3 @ 3.50 GHz core with 32 GB of memory.

In the analytical and aircraft optimization examples in the following sections, we use the following notation: 
\begin{itemize}
    \item[a)]  Optimization algorithms:
\begin{itemize}
    \item \texttt{VBO}: a vanilla Bayesian optimization per~\Cref{EGOalgorithm}. 
    \item \texttt{TLBO}: Bayesian optimization with transfer learning per~\Cref{TLBO-algo}.
\end{itemize}
 \item[b)] Surrogate models within \texttt{TLBO}:
\begin{itemize}
    \item \texttt{ETL}: ensemble of surrogates using transfer learning per~\Cref{TL-ensemble-algo}.  
\end{itemize}
 \item[c)] Variance estimation for ensemble methods:
\begin{itemize}
    \item \texttt{TV}: Target variance using from~\Cref{eq:variance-target}.
    \item \texttt{AV}: Weighted variance based on ensemble of surrogates without clustering using~\Cref{eq:variance-ensemble}.
\end{itemize}
\end{itemize}
Combinations of the above methods are concatenated when presented. For example, \texttt{TLBO} optimizer with \texttt{ETL} and \texttt{AV} is presented as \texttt{TLBO-ETL-AV}.

We use a single cluster assumption for creating the ensembles of surrogates in \texttt{ETL} to minimize its impact on the computational costs of optimizations. We set fixed weights for the ensemble criteria from~\Cref{sec:multi-criteria}. For the objective function, we are interested in the shape of the source functions and how they compare to the target function. Therefore, we set the criteria weights to be based on shape and variance. Whereas for the constraint functions, accuracy is just as important as the shape of the function as it determines feasibility in the acquisition function maximization loop as described in~\Cref{sec:constraint-handling}. Therefore, we assign equal weights to each of the shape, accuracy, and variance criteria.
\subsection{Analytical examples}
We first illustrate an unconstrained optimization problem based on the Bohachevsky functions~\cite{bohachevsky1986generalized, DBLP:journals/corr/JamilY13}. We use three homogeneous source functions, $f_{\mbox{S1}}$, $f_{\mbox{S2}}$, $f_{\mbox{S3}}$, and one target function, $f_{\mbox{T}}$, as follows. For a given $x_1, x_2 \in [-5, 5]$, we consider
\begin{align}
f_{\mbox{S1}}(x_1, x_2) &= x_1^2 + 2x_2^2 - 0.3\cos(\pi x_1) - 0.4\cos(2\pi x_2) + 0.7 \label{eq:tl-f1}\\
f_{\mbox{S2}}(x_1, x_2) &= x_1^2 + 2x_2^2 - 0.3\cos(\pi x_1)\cos(2\pi x_2) \label{eq:tl-f2}\\
f_{\mbox{S3}}(x_1, x_2) &= 2x_1^2 + 4x_2^2 - 0.3\cos(3\pi x_1 + 4\pi x_2) - 0.5 \label{eq:tl-f3}\\
f_{\mbox{T}}(x_1, x_2) &= 0.5x_1^2 + x_2^2 - 0.3\cos(3\pi x_1 + 4\pi x_2) + 0.4 .\label{eq:tl-f4}
\end{align}

\begin{figure}[!h]
    \centering
    \begin{subfigure}{.45\textwidth}
    \centering 
    \includegraphics[width=6 cm]{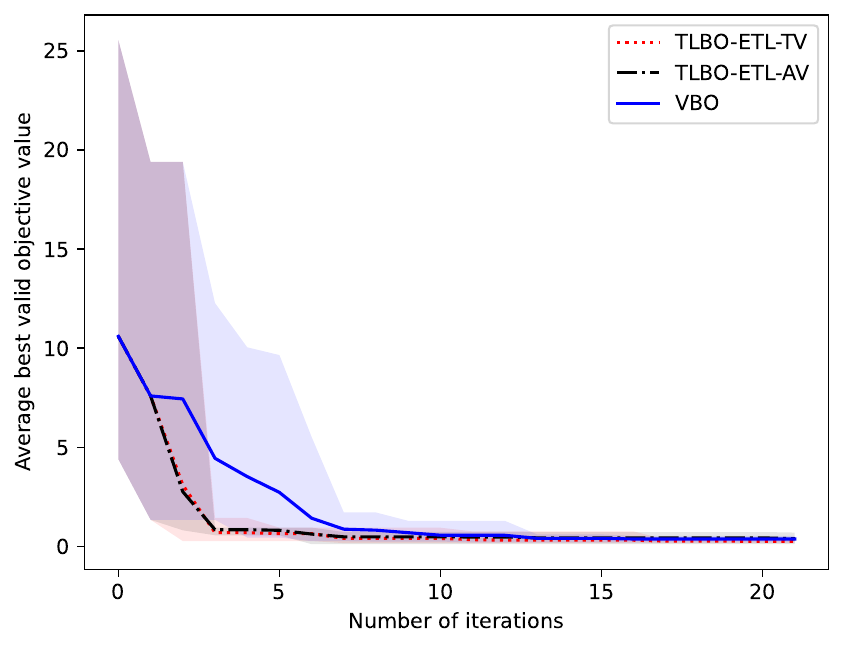}
    \caption{Convergence plots}   
    \label{TLBO-analytical1-convergnce-figurea}
    \end{subfigure}
    \hfill
    \begin{subfigure}{.45\textwidth}
    \centering
	\includegraphics[width=6 cm]{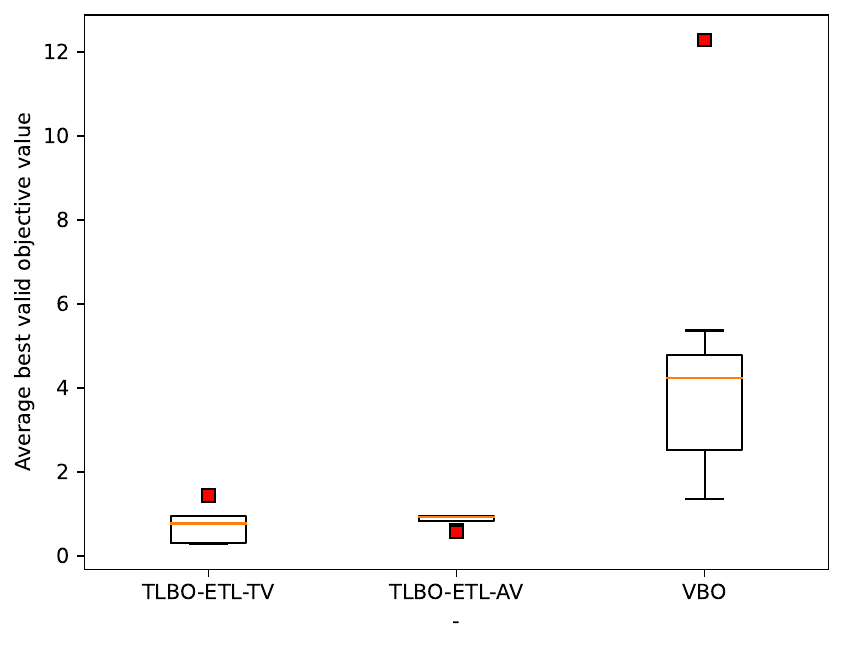}
    \caption{Box plots at iteration 5}
    \label{TLBO-analytical1-convergnce-figureb}
    \end{subfigure}
    \caption{Results of the analytical problem showing the improved convergence using \texttt{TLBO} over the standard \texttt{VBO} over 20 optimization runs.}
    \label{TLBO-analytical1-convergnce-figure}
\end{figure}

We conduct a DOE of size 50 using Latin hypercube sampling (LHS) for each of the source problems to create the source surrogate models. Then, we initiate the Bayesian optimization algorithm using an initial DOE of the target function of 2 random evaluations. We set the number of iterations to 20 and perform optimizations using both \texttt{VBO} and \texttt{TLBO} using \texttt{ETL}. For \texttt{TLBO} we also test two possible methods for variance calculation: \texttt{TV} and \texttt{AV}. We repeat the process for 20 runs to obtain the convergence results presented in~\Cref{TLBO-analytical1-convergnce-figure}. We note that \texttt{TLBO} using \texttt{ETL} achieves faster convergence than \texttt{VBO} as shown in ~\Cref{TLBO-analytical1-convergnce-figurea}. We note that both \texttt{TV} and \texttt{AV} yield similar convergence results, as seen in the box plots in~\Cref{TLBO-analytical1-convergnce-figureb}, which can be explained by the fact that all source surrogate models are built based on LHS sampling across the same design space as the target problem. For the remainder of this work, we use \texttt{TV} to maintain acquisition function exploration based on the target model. 

We also tested our methodology using a benchmark problem typically used for multi-fidelity optimization~\cite{mainini2022analytical}. The problem is based on variations of the Rosenbrock function $h_{\mbox{T}}$ by creating a medium fidelity level $h_{\mbox{S1}}$ and a low fidelity level $h_{\mbox{S2}}$. In this transfer learning scenario, we consider the low and medium fidelity functions as source problems and the higher fidelity function as the target problem. 
\begin{eqnarray}
     h_{\mbox{S1}}(x) &=&  \sum _ {i=1}^ {D-1}  50  (x_ {i+1}-x_ {i}^ {2})^ {2}  +  (-2-x_ {i})^ {2}  -  \sum _ {i=1}^ {D}  0.  5x_ {i} \label{rosen:1}\\ 
       h_{\mbox{S2}}  (x)&=&  \frac {f_ {1}(x)-4-\sum _ {i=1}^ {D}0.5x_ {i}}{10+\sum _ {i=1}^ {D}0.25x_ {1}} \label{rosen:2} \\
         h_{\mbox{T}}  (x)&=&  \sum _ {i=1}^ {D-1}  100  (x_ {i+1}-x_ {i}^ {2})^ {2}  +  (1-x_ {i})^ {2}. \label{rosen:3}
\end{eqnarray}

\begin{figure}[h!]
	
    \begin{subfigure}{0.45\textwidth}
    \hspace*{-1in}
    \centering 
    \includegraphics[width=6.5 cm]{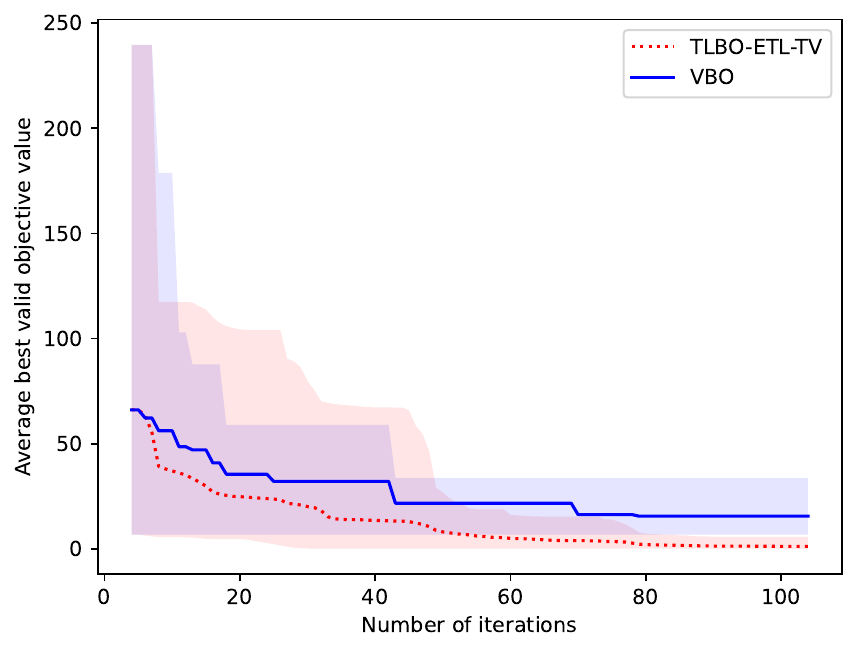}
    \caption{Convergence plots relative to the number of iterations}   
    \end{subfigure}
    \hspace*{0.2in}
    \begin{subfigure}{0.45\textwidth}
    \centering
	\includegraphics[width=6.5 cm]{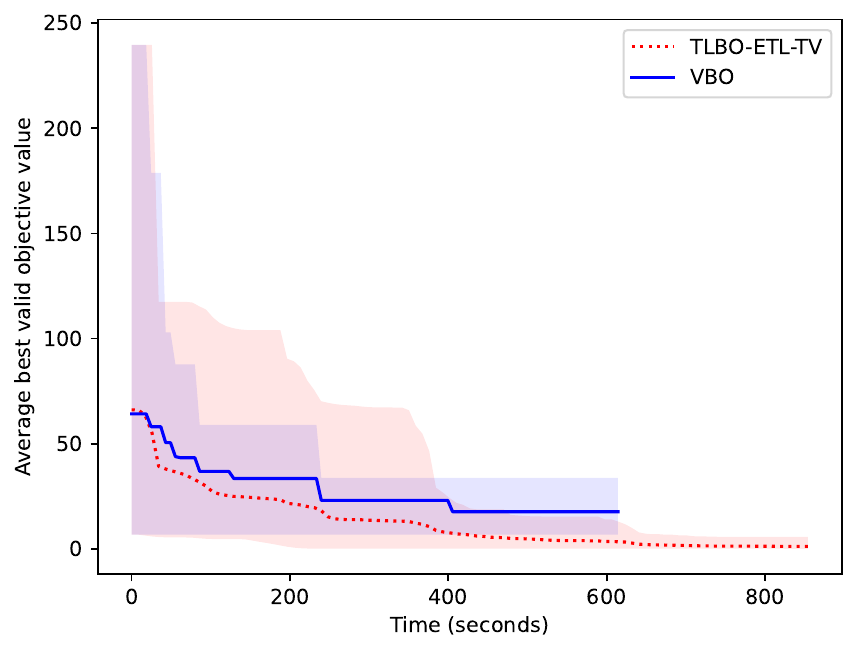}
    \caption{Convergence plots relative to computational time}
    \label{fig:rosen-comp-cost}
    \end{subfigure}
    \begin{subfigure}[t]{1\textwidth}
   \centering
	\includegraphics[width=6.5 cm]{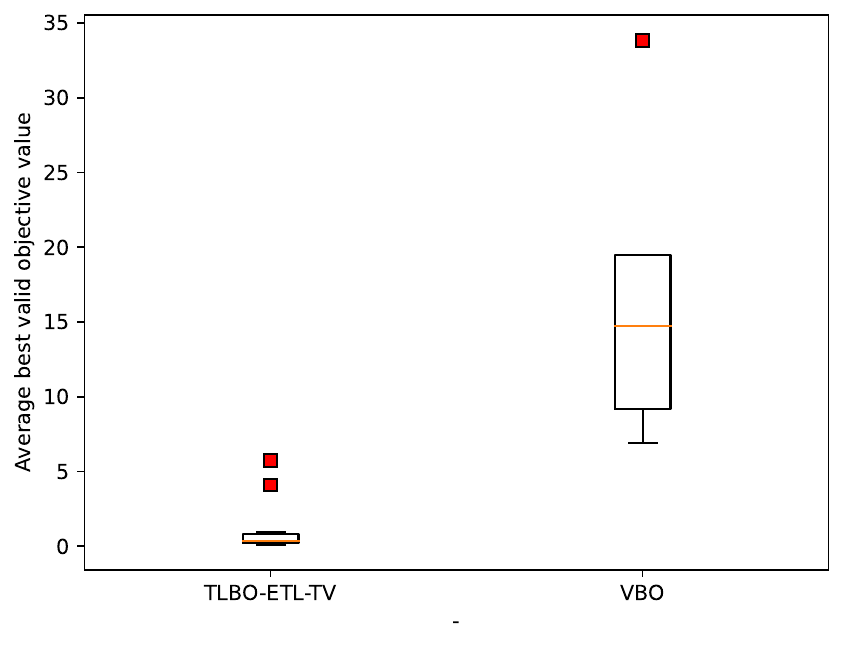}
    \caption{Box plots of the results}
    \label{fig:rosen-comp-boxplot}
    \end{subfigure}
    
    \caption{Results of the 5-dimensional Rosenbrock multi-fidelity problem  showing the improved convergence of \texttt{TLBO} over \texttt{VBO}  over 20 optimizations runs.}
    \label{TLBO-rosenbrock5d-convergnce-figure}
\end{figure}

For this problem, we select $D=5$ and $-2\leq x_i \le 2$ for $i = 1,...,D$ as the so-called \texttt{MF2.2} benchmark problem from~\cite{mainini2022analytical}. We create a DOE of size 100 using LHS  using each of $h_{\mbox{S1}}$ and $h_{\mbox{S2}}$. We set the size of the initial DOE for the target problem as the number of design variables $d$. We conduct 20 optimization runs of 100 iterations each and compare \texttt{TLBO} and \texttt{VBO}. We present the results in~\Cref{TLBO-rosenbrock5d-convergnce-figure} showing the improved convergence when using \texttt{TLBO} over \texttt{VBO}. In~\Cref{fig:rosen-comp-cost}, we assume that the computational cost of each high fidelity function evaluation is 10 seconds. It is noted that \texttt{TLBO} produces lower average minima throughout the optimization even with the added computational costs as seen in~\Cref{fig:rosen-comp-boxplot}. 

\subsection{Aircraft conceptual design optimization}
\label{sec:tlbo-aircraft}
We use an aircraft conceptual design problem to demonstrate our Bayesian optimization transfer learning approach. We present a scenario based on the problem defined in~\Cref{sec:into-aircraft-problem}
We select a low number of evaluations for the initial DOE of the target optimization (less than a quarter of the dimension of the input vector $x$), and 150 evaluations for the DOEs of the source problems. Note that DOE size does not have to be equal across all source problems. 

\begin{figure}[H]
    \centering
	\begin{subfigure}{.45\textwidth}
    \centering 
    \hspace*{-1in}
    \includegraphics[width=6.5 cm]{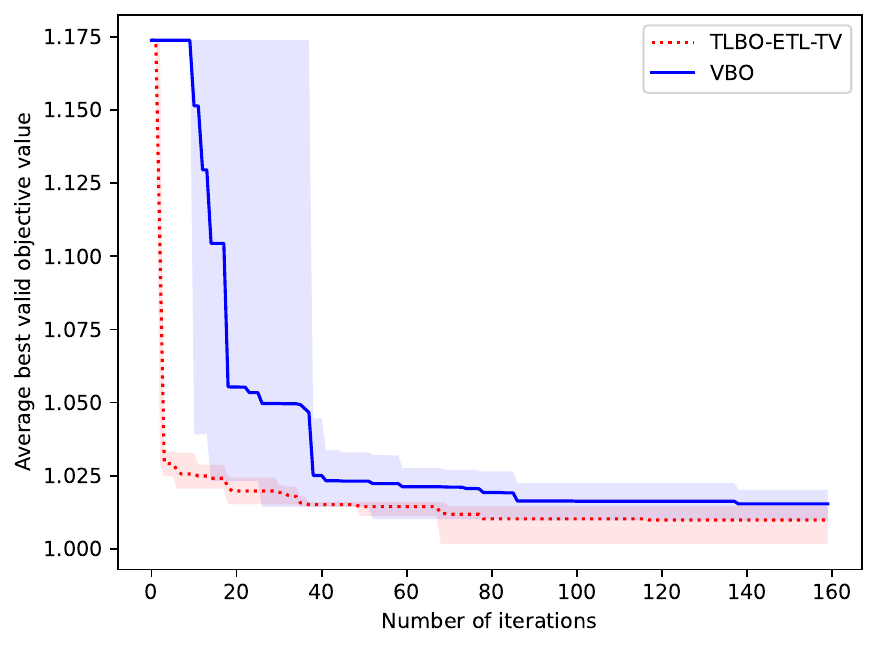}
    \caption{Convergence plots relative to number of iterations}   
    \label{tlbo-ac-conv-iter1}
    \end{subfigure}
    \hspace*{0.2in}
    \begin{subfigure}{.45\textwidth}  
    \centering
    \includegraphics[width=6.5 cm]{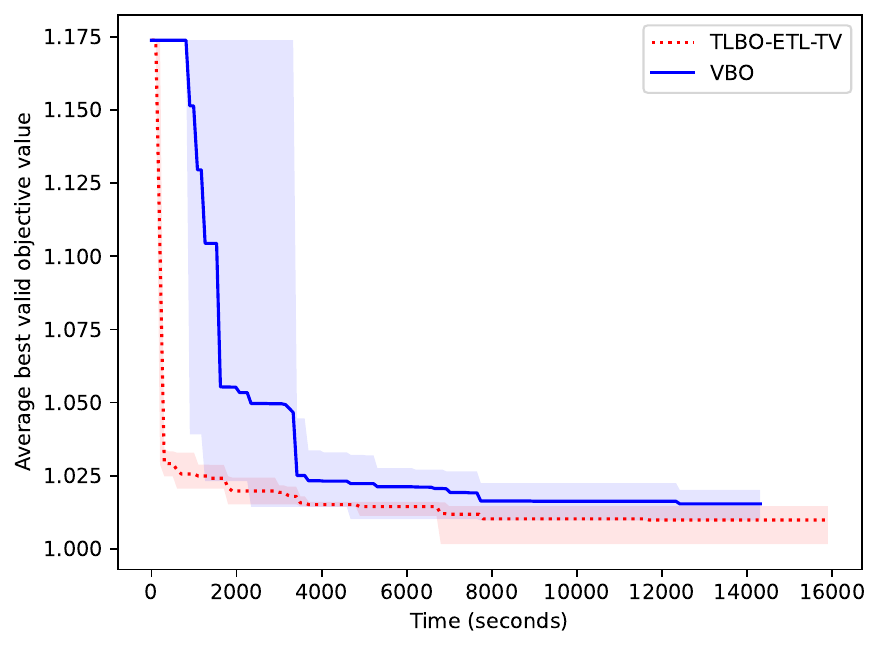}
    \caption{Convergence plots relative to computational time}
    \label{tlbo-ac-conv-time}
    \end{subfigure}
 \begin{subfigure}{1\textwidth}
    \centering
	\includegraphics[width=6.5 cm]{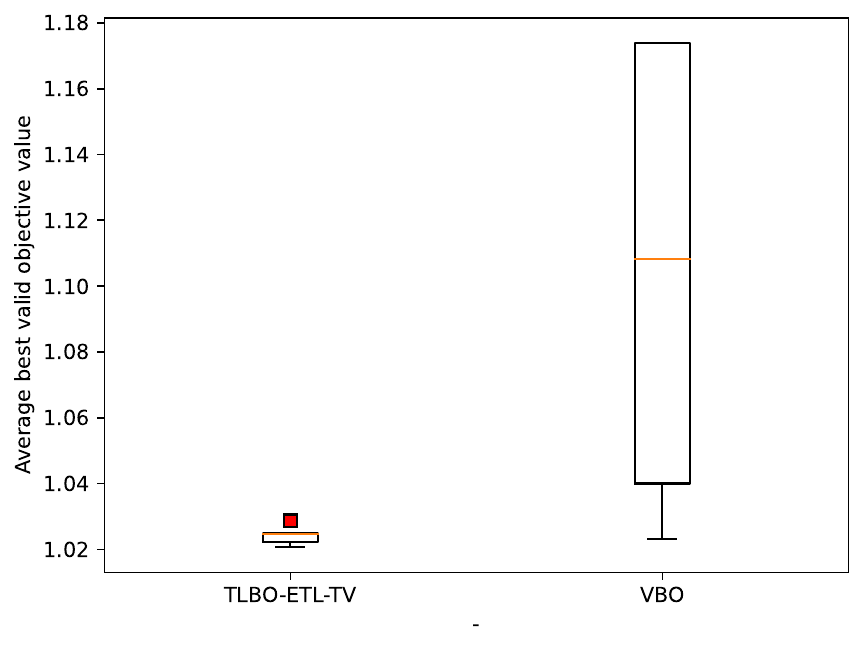}
    \caption{Box plots of the best valid minima at iteration 15}
    \label{tlbo-ac-boxplot}
    \end{subfigure}
    \caption{Results of 10 runs of the BRAC optimization problem showing the improved convergence and the variance of \texttt{TLBO} over \texttt{VBO}.}
    \label{TLBO-ac-results1-figure}
\end{figure}

We first perform 10 sets of the initial DOEs consisting of 3 samples each using the target BRAC optimization problem. We use these DOEs to conduct 10 optimization runs for each of \texttt{TLBO} and \texttt{VBO}. We present the results of the optimizations in~\Cref{TLBO-ac-results1-figure} demonstrating the superior performance of \texttt{TLBO} over \texttt{VBO}. Convergence plots show the average of the minimum normalized objective values at every iteration (shown as a line per method) and the variance at every iteration (shown as pastel colors of the corresponding line). We note the improved convergence behaviour of \texttt{TLBO} in the early iterations of the optimization even when considering the additional computation cost as can be seen in~\Cref{tlbo-ac-conv-time}. This improvement is also demonstrated in~\Cref{tlbo-ac-boxplot} where we show result boxplots of the two tested methods at the 15\textsuperscript{th} optimization iteration. 

\begin{figure}[H]
    \centering
    \captionsetup[subfigure]{justification=centering}
	\hspace{-.3cm}
    \begin{subfigure}[t]{0.32\textwidth}
    \centering 
    \includegraphics[width=5 cm]{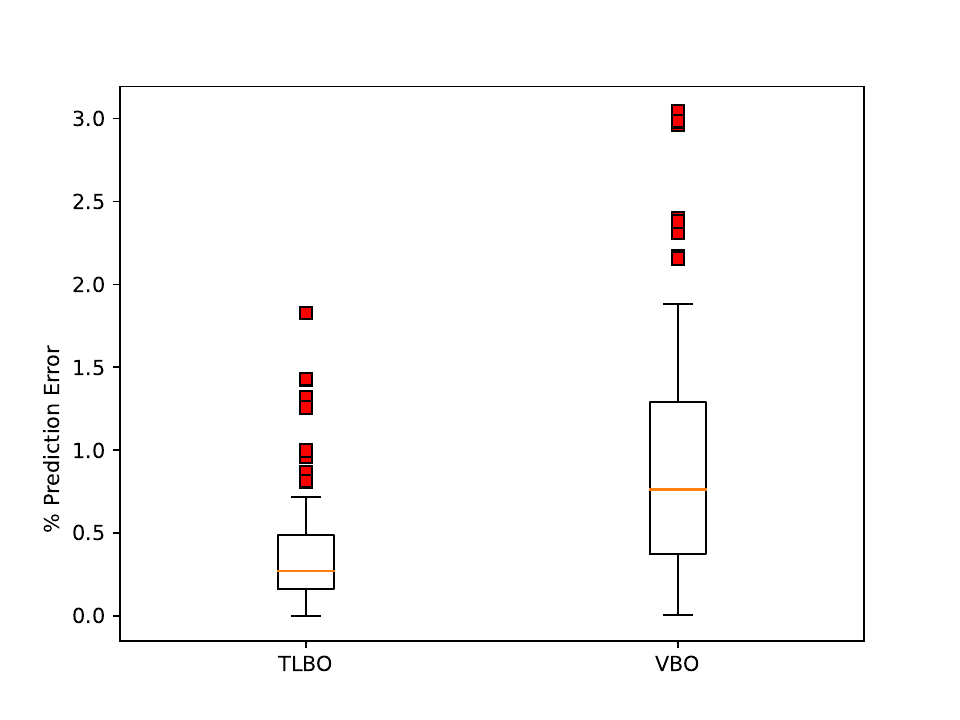}
    \caption{MTOW (objective) }   
    \end{subfigure}
~~\begin{subfigure}[t]{0.32\textwidth}
    \includegraphics[width=5cm]{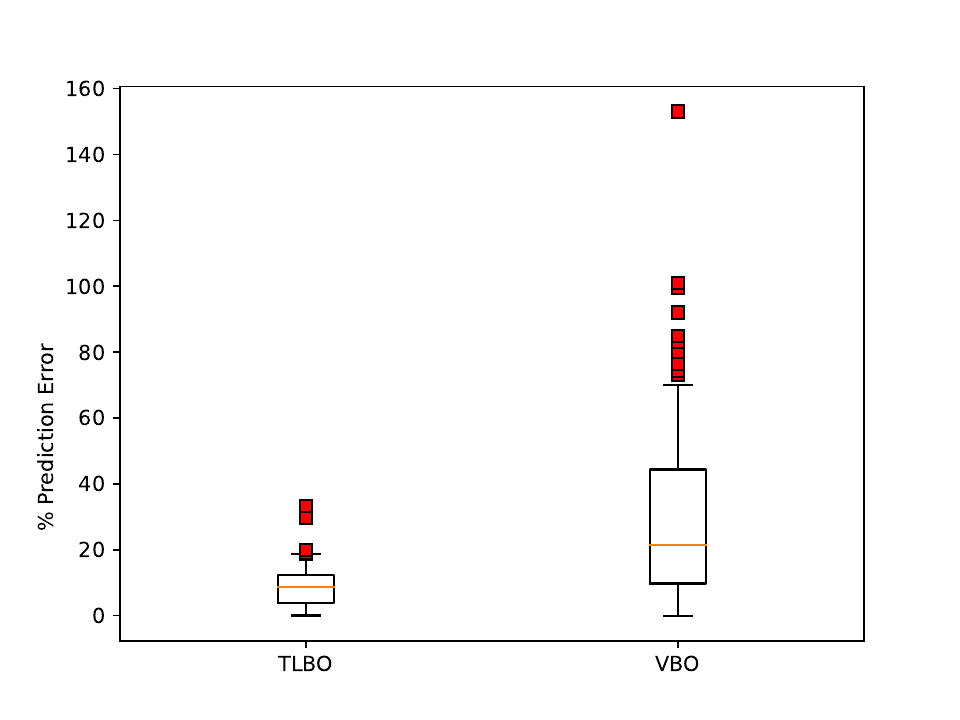}
    \caption{ $c_1$ constraint\textsuperscript{*}(BFL)}
    \end{subfigure}
~~\begin{subfigure}[t]{0.32\textwidth}
    \centering 
    \includegraphics[width=5 cm]{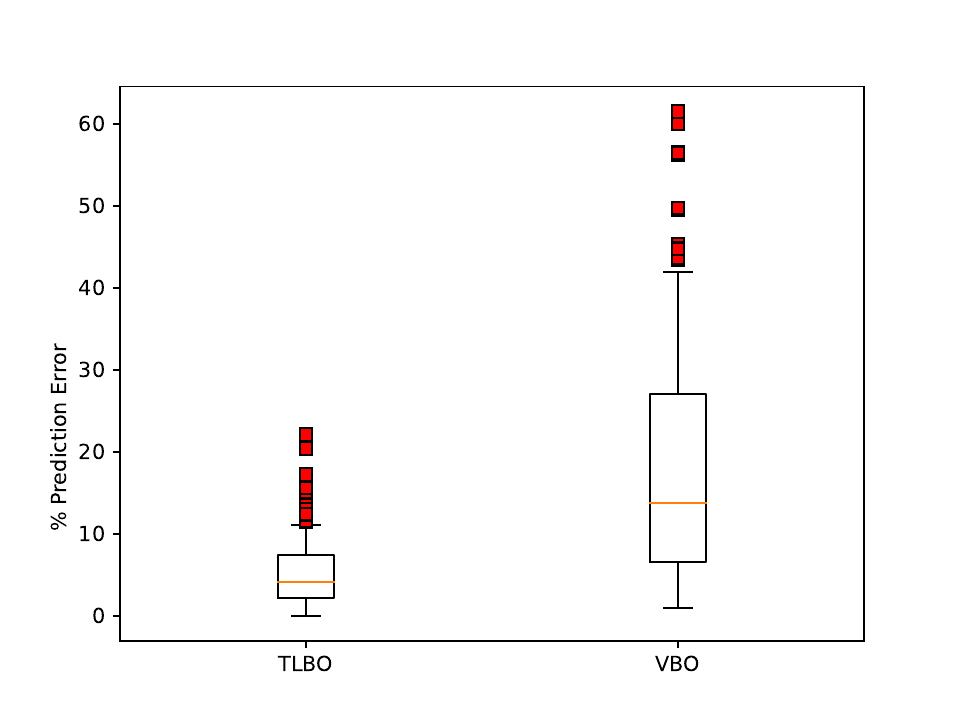}
    \caption{ $c_2$ constraint (ICA)}  
    \end{subfigure}
   
  	\hspace{-.3cm}
  \begin{subfigure}[t]{0.32\textwidth}
    \centering
    \includegraphics[width=5 cm]{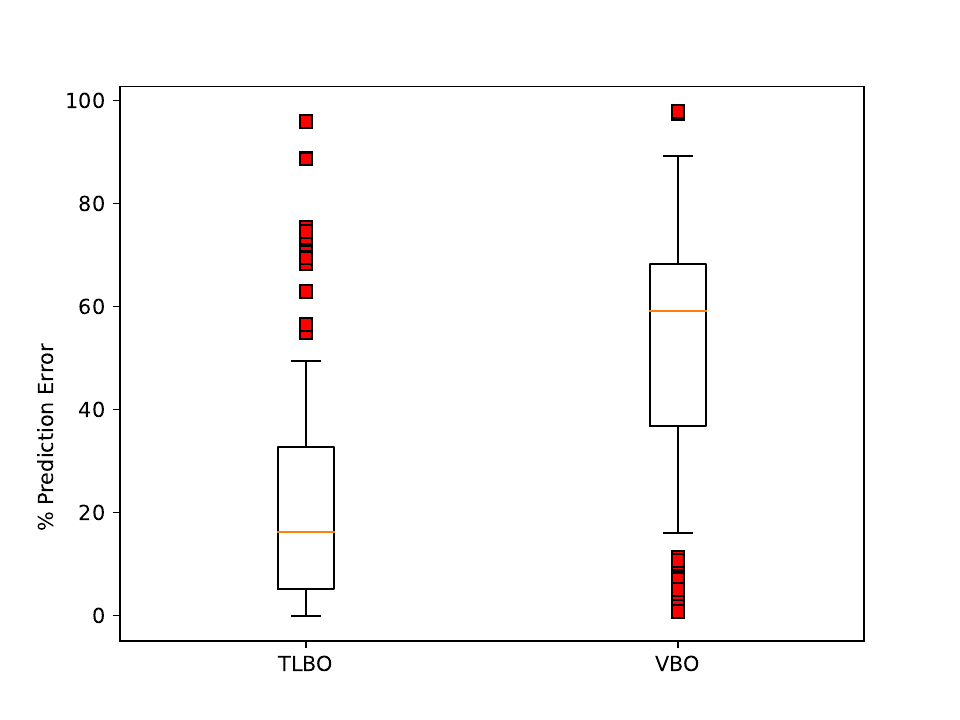}
    \caption{ $c_3$ constraint ($V_{\mbox{ref}}$)}
    \end{subfigure}
    ~~\begin{subfigure}[t]{0.32\textwidth}
    \centering 
    \includegraphics[width=5 cm]{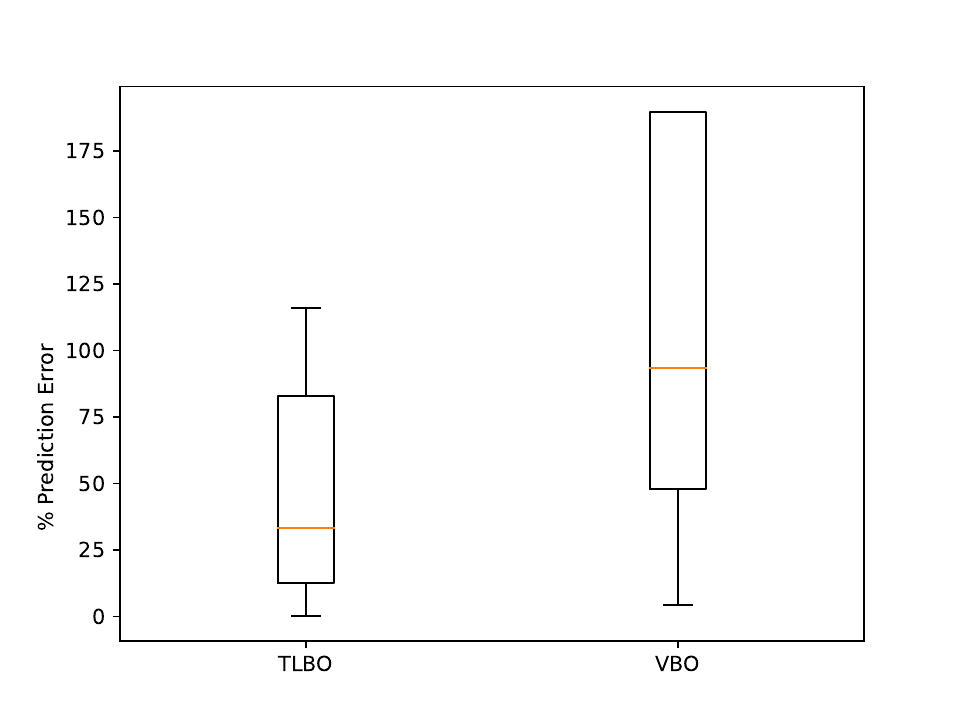}
    \caption{ $c_5-c_8$ constraints (Flight controls)}
    \end{subfigure}  
    ~~\begin{subfigure}[t]{0.32\textwidth}
    \includegraphics[width=5 cm]{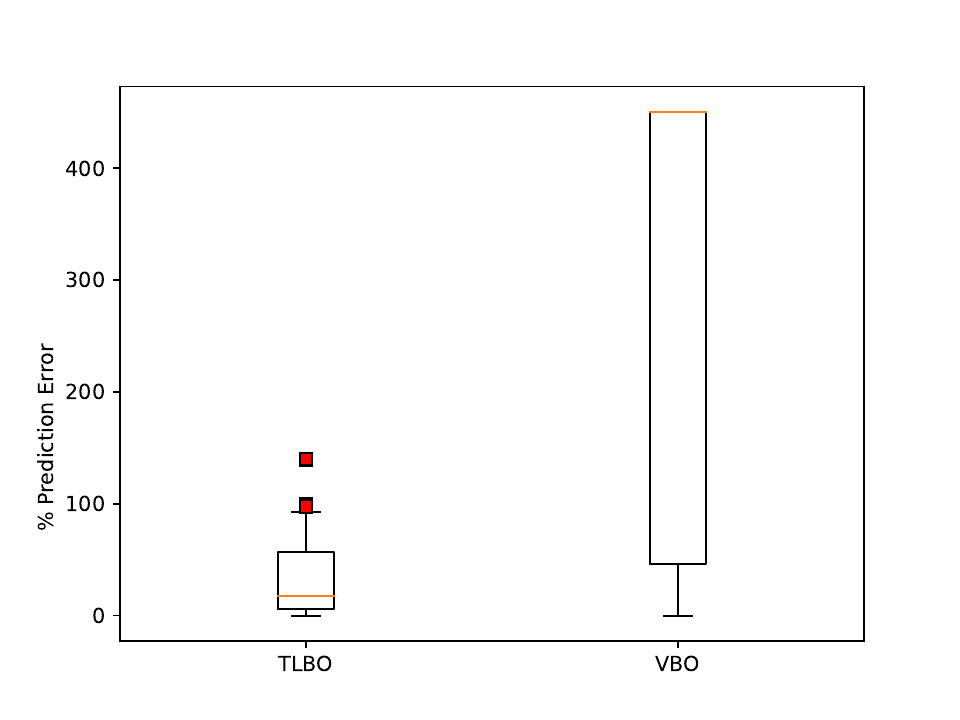}
    \caption{$c_9$ constraint\textsuperscript{*} (Landing gear )}
    \label{TLBO-box-plot-ac-1}       
    \end{subfigure}
    
    	\hspace{-.3cm}\begin{subfigure}[t]{0.32\textwidth}
    
    \includegraphics[width=5cm]{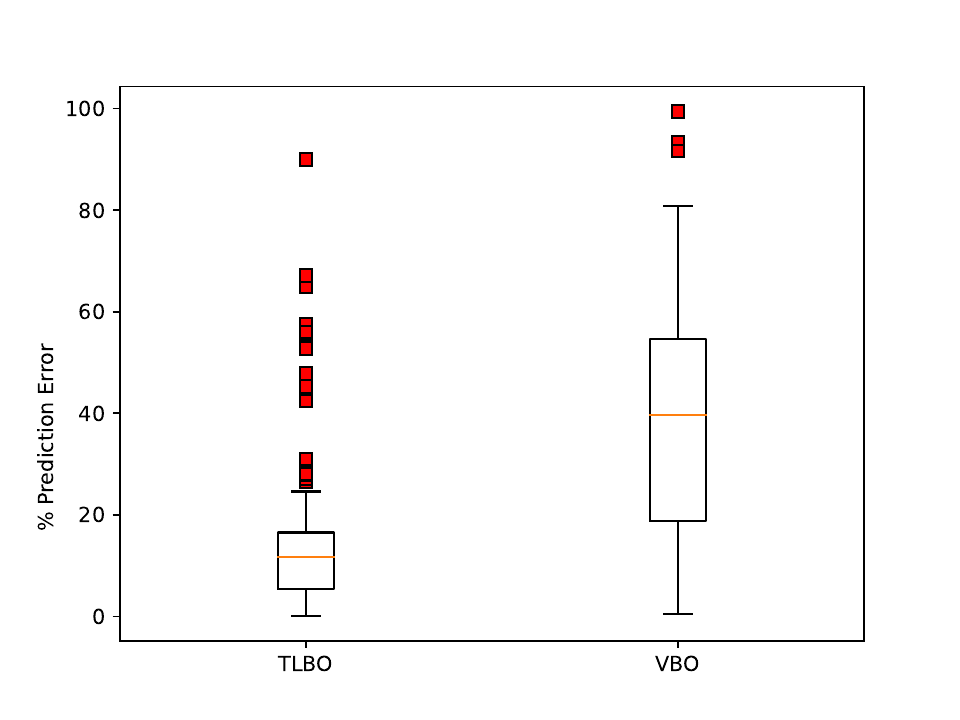}
    \caption{$c_4$ constraint (Excess fuel) }   
    \end{subfigure}   ~~
    \begin{subfigure}[t]{0.32\textwidth}
    \centering 
    \includegraphics[width=5 cm]{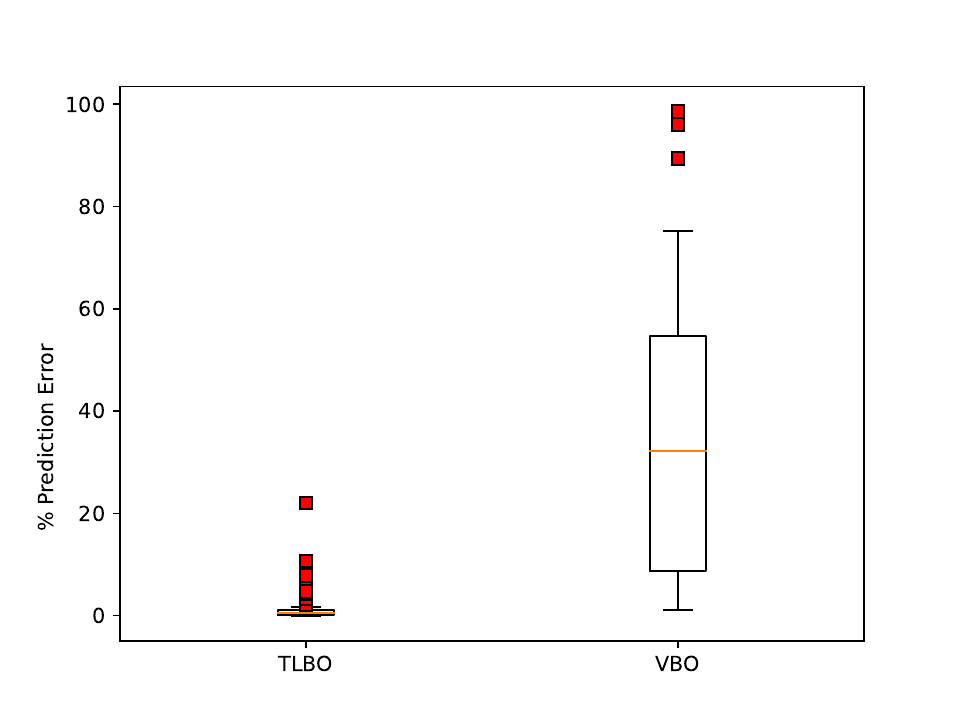}
    \caption{$c_{12}$ constraint\textsuperscript{*} (Range)}

    \end{subfigure}
    ~~
\begin{subfigure}[t]{.32\textwidth}
  
    \includegraphics[width=5 cm]{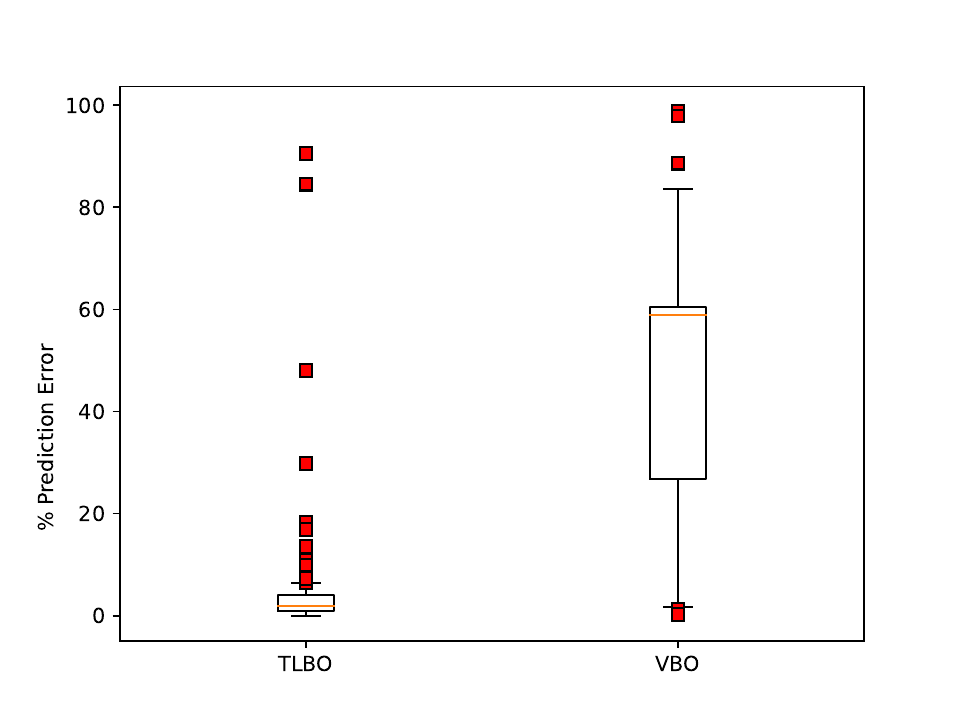}
    \caption{$c_{10}$ constraint\textsuperscript{*} (Tip chord)}  
    \end{subfigure}    
     \hspace{-0.0in}\raggedright
    
    \caption{Box plots comparing the prediction errors relative to the objective function and the constraints of \texttt{TLBO} and \texttt{VBO} at iteration 15 of the optimization. }
    \footnotesize{\textsuperscript{*} denotes active inequality constraints (less than 2\% margin) at the minimum observed valid objective function  using 100 test evaluations.}
    \label{TLBO-analytical1-convergnce-figure1}
\end{figure}

We compare the prediction accuracy of the surrogate models in the Bayesian optimization frameworks at the 15th\textsuperscript{th} iteration using 100 test points in~\Cref{TLBO-analytical1-convergnce-figure1}, showing the improved accuracy of \texttt{TLBO} across all the surrogate models for constraints and the objective function, which explains the improved convergence results early in the optimization per~\Cref{tlbo-ac-conv-iter1}.

\subsection{Discussion}
While the proposed \texttt{TLBO} method demonstrates improved convergence in the presented examples, practitioners should be aware of scenarios where the method may underperform or fail to provide benefits over standard Bayesian optimization.

\paragraph{Highly dissimilar source and target problems.} When source data originates from fundamentally different optimization landscapes—for example, when the source problem is unimodal while the target is highly multimodal—the transferred surrogate models may mislead the optimization search. Although the multi-criteria weighting scheme assigns lower probabilities to poorly matching sources, sufficiently dissimilar problems may still introduce bias in the ensemble predictions, particularly in early iterations when limited target data is available for weight calibration.

\paragraph{Conflicting source models.} When multiple source models disagree significantly and suggest optima in different regions of the design space, the weighted ensemble produces an averaged prediction that may not accurately represent any actual optimum. This averaging effect can result in a surrogate landscape that is smoother than the true target function, potentially causing the optimizer to miss local optima or converge to suboptimal regions.

\paragraph{Insufficient target data for weight calibration.} The multi-criteria weighting relies on the target DOE to assess source model quality. With very few initial target samples (e.g., fewer than the number of design variables), the weight estimation may be unreliable, leading to inappropriate source model rankings. In such cases, the alternating \texttt{TLBO}/\texttt{VBO} procedure described in Section~\ref{sec:negative-transfer} can help mitigate this risk.

\paragraph{Constraint boundary shifts.} When constraints are active in different regions between source and target problems---for instance, if a structural constraint becomes critical at different wing loading conditions---transferring the constraint surrogate may incorrectly predict feasibility boundaries. This is particularly problematic for inequality constraints near their bounds, where small prediction errors can lead to misclassification of feasible and infeasible regions.

\paragraph{Computational overhead.} The ensemble of surrogates approach adds computational cost due to source model construction, weight calculation, and ensemble predictions at every iteration. For simple low-dimensional problems, fast blackbox evaluations, or when source data quality is poor, this overhead may not be justified by convergence improvements. Users should assess whether the expected reduction in blackbox evaluations outweighs the added surrogate modeling costs.

\paragraph{Addressing Negative Transfer}
\label{sec:negative-transfer}
Negative transfer occurs when source data degrades target optimization performance, typically when source and target problems are insufficiently related. Our methodology mitigates negative transfer through several mechanisms. First, the multi-criteria weighting scheme (\Cref{sec:multi-criteria}) automatically assigns lower probabilities to source models that poorly match the target data, effectively down-weighting irrelevant sources. Second, the scaling and bias adjustment in the transfer of prior (\Cref{TL-prior-alphab}) allows source models to adapt to the target function's range, reducing the impact of value mismatches. Third, for cases where negative transfer may still occur, an alternating procedure can be employed: the optimization alternates between \texttt{TLBO} and the standard Bayesian optimization (\texttt{VBO}) at specified intervals, allowing the algorithm to rely on target-only surrogates when the transfer learning ensemble under performs. This hybrid approach ensures robustness by periodically resetting to vanilla Bayesian optimization, thereby bounding the potential negative impact of misleading source information while still benefiting from positive transfer when source data is informative. 
\newpage
\section{Conclusions}
\label{sec:TLBO-conclusions}
In this work, we advanced the state-of-the-art for transfer learning in Bayesian optimization by presenting a novel method to address constraints and heterogeneity. Our method is based on an ensemble of surrogates using transfer learning algorithm which we embedded in a Bayesian optimization framework. We propose two approaches to handle heterogeneity between source and target optimization problems: the first is based on meta data adjusted dimension reduction aimed to solve design space heterogeneity, and the second approach addresses heterogeneity of constraints by selecting which constraints are to be used as source constraints in a target optimization problem. The method adds computational costs to the standard Bayesian optimization algorithm due to the additional source model selection step at every iteration. We validated the proposed methodology using an aircraft conceptual design optimization scenario and compared the results with standard Bayesian optimization. The results show significant improvement in convergence early in the optimization even with the increased computational cost due to the introduction of the ensemble of surrogates. 

Future work includes extending the weighted ensemble to a mixture of experts approach with multiple clusters, where each cluster can have different source model weights based on local performance as presented for surrogate modeling in~\cite{tfaily2026transfer}. This would allow the ensemble to adapt spatially across the design space, potentially improving performance in regions where different source models are more accurate.

Additional future research directions include: extending the application to more aircraft design optimization test cases with varying fidelity of simulation models to understand the impact of source data fidelity on optimization results; investigating adaptive schemes that automatically select between \texttt{TLBO} and standard Bayesian optimization based on detected negative transfer; and developing methods to quantify source-target similarity \textit{a priori} to guide practitioners in selecting appropriate source data for transfer. 

\backmatter

\bmhead{Acknowledgements}

The authors would like to thank Jasveer Singh from Bombardier for his insights and support in setting up the industrial aircraft conceptual design problems.  

\section*{Declarations}
\bmhead{Author contributions}
All authors contributed to the proposed methods and presented results. The first draft of the manuscript was written by the first author and all authors read, improved, and approved the submitted manuscript.

\bmhead{Funding}
The first and last authors are grateful for the partial support of this work by NSERC grant (RGPIN 436193-24);  this support does not constitute an endorsement of the opinions expressed in this paper. Funding for the third author's research 
was partially provided by the NSERC Discovery grant (RGPIN-2024-0509).

\bmhead{Data availability}
Due to Bombardier intellectual property considerations, the aircraft application case data cannot be shared.

\bmhead{Competing interests}
The authors declare that there is no conflict of interest.

\bmhead{Replication of results}
This paper has provided all the necessary information and equations to reproduce the analytical results. Due to Bombardier intellectual property considerations, the simulation software to produce industrial application case results cannot be shared. 

\bmhead{Ethics approval}
Not applicable

\bmhead{Consent to participate}
Not applicable









\bibliography{references}

\end{document}